\documentclass[11pt,a4paper]{article}
\usepackage[english]{babel}
\usepackage{eurosym}
\usepackage{amsfonts}
\usepackage{amssymb}
\usepackage{epsfig}
\usepackage{graphicx}
\usepackage{amssymb,amsfonts,amsmath,enumitem,mathrsfs,amsthm,mathtools}
\makeatletter
\newcommand{\leqnomode}{\tagsleft@true\let\veqno\@@leqno}
\newcommand{\reqnomode}{\tagsleft@false\let\veqno\@@eqno}
\makeatother

\usepackage{epstopdf} 
\usepackage{tikz}
\usetikzlibrary{intersections}
\usetikzlibrary{arrows,arrows.meta,bending}

\usepackage{cite}

\usepackage{xcolor}
\definecolor{candyapplered}{rgb}{1.0, 0.03, 0.0}
\definecolor{coolblack}{rgb}{0.0, 0.18, 0.39}
\definecolor{darkpowderblue}{rgb}{0.0, 0.2, 0.6}
\definecolor{lightbrown}{rgb}{0.71, 0.4, 0.11}
\definecolor{ruddybrown}{rgb}{0.73, 0.4, 0.16}
\usepackage{hyperref}
\hypersetup{hyperfigures = true, colorlinks = true, linkcolor=darkpowderblue, citecolor=ruddybrown}

\usepackage[utf8]{inputenc}

\def\Sum{\displaystyle\sum}

\def\Inf{\displaystyle\inf}

\newcommand{\pr}{\medskip\noindent\textit{\textbf{Proof.}} }

\newcommand{\prrr}[1]{\medskip\noindent\textit{\textbf{Proof of Theorem #1.}} }

\newtheorem{theo}{Theorem}[section]
\newtheorem{prop}{Proposition}[section]
\newtheorem{lem}{Lemma}[section]
\newtheorem{rem}{Remark}[section]
\newtheorem{defi}{Definition}[section]

\begin{document}
	\title{Existence results for singular elliptic problem involving a fractional $p$-Laplacian.}
	\author{Hana\^{a} ACHOUR$^{a}$ and Sabri BENSID$^{b}$ }
	\maketitle
	\begin{center}
		$^{a}$Dynamical Systems and Applications Laboratory\\
		Department of Mathematics, Faculty of Sciences\\
		University of Tlemcen, B.P. 119\\
		Tlemcen 13000, Algeria\\
		Mail: $hanaa495@outlook.com$
		
		$^{b}$Dynamical Systems and Applications Laboratory\\
		Department of Mathematics, Faculty of Sciences\\
		University of Tlemcen, B.P. 119\\
		Tlemcen 13000, Algeria\\
		Mail: $ edp\_sabri@yahoo.fr$
	\end{center}
	\addcontentsline{toc}{section}{Abstract}
	\begin{abstract}
		In this article, the problems to be studied are the following
		\leqnomode
		\begin{equation*}
			\label{p}
			\left\{\begin{array}{ll} (-\Delta )_p^s u \pm \dfrac{|u|^{p-2}u}{|x|^{sp}} = \lambda f(x,u) & \quad \mbox{in }\ \Omega\\[0.3cm]
				u= 0        & \quad \mbox{on }\ \mathbb{R}^N \setminus \Omega,\tag{P$_{\pm}$}
			\end{array}
			\right.
		\end{equation*}
		\reqnomode
		where $\Omega$ is a bounded regular domain in $\mathbb{R}^N(N\geq 2)$ containing the origin, $p>1$, $s\in(0,1)$, $(N>ps)$, $\lambda>0$,
		$f : \Omega \times \mathbb{R} \longrightarrow \mathbb{R}$ is a Carathéodory function satisfying a suitable growth condition and $(-\Delta )_p^s$ is the fractional p-Laplacian defined as
		$$(-\Delta )_{p}^{s} u(x) \coloneqq  \displaystyle 2 \lim_{\varepsilon \rightarrow 0} \int_{\mathbb{R}^N \setminus B_{\varepsilon}(x)} \dfrac{\vert u(x)-u(y) \vert^{p-2}(u(x)-u(y))}{\vert x-y \vert^{N+sp}} ~dy, ~~~~ x \in \mathbb{R}^N,$$
		where $B_{\varepsilon}(x)$ is the open $\varepsilon$-ball of centre $x$ and radius $\varepsilon$. Using the critical point theory combining to the fractional Hardy inequality, we show that the problem $(P_+)$ admits at least two distinct nontrivial weak solutions. For the problem $(P_-),$ we use the concentration-compactness principle for fractional Sobolev spaces to give a weak lower semicontinuity result and prove that problem $(P_-)$ admits at least one non-trivial weak solution.
		
	\end{abstract}
	\newline
	\newline
	\noindent {\textbf{Keywords:}}{ Singular problem, fractional $p$-Laplacian, critical point, variational method, fractional Hardy inequality.}\hfill \break
	
	\section{Introduction}
	Singular elliptic problems arise in several context like chemical heterogeneous catalysts, in the theory of pseudoplastic fluid, non-Newtonian fluid and in the study of relativistic matter in magnetic fluid. for more details, see \cite{Nachman1980,Diaz1987,Diaz1985}.
	
	On the other hand, recently, a great deal of works is devoted to the fractional nonlocal problem arising in many fields such as anomalous diffusion in plasma, flames propagation, geophysical fluid dynamics and American options in finances. See \cite{Cabre2010,Dinezza2012} and the references therein.
	
	In this work, we continue the study of semilinear singular elliptic problem with the nonlocal operator theory. more precisely, we are interested by the fractional $p$-Laplacian operator which up to normalization functions may be  defined as
	$$(-\Delta )_{p}^{s} u(x) \coloneqq  \displaystyle 2 \lim_{\varepsilon \rightarrow 0} \int_{\mathbb{R}^N \setminus B_{\varepsilon}(x)} \dfrac{\vert u(x)-u(y) \vert^{p-2}(u(x)-u(y))}{\vert x-y \vert^{N+sp}} ~dy, ~~~~ x \in \mathbb{R}^N,$$
	where $B_{\varepsilon}(x)$ is the open $\varepsilon$-ball of centre $x$ and radius $\varepsilon$. We refer to \cite{Caffarelli2012,Iannizzotto2016} for more details.
	
	In the last years, many works are devoted to the quasi-linear problem
	\begin{equation}
		\label{q}
		\left\{\begin{array}{ll} (-\Delta )_p^s u = f(x,u) & \quad \mbox{in }\ \Omega\\[0.3cm]
			u= 0        & \quad \mbox{on }\ \mathbb{R}^N \setminus \Omega.
		\end{array}
		\right.
	\end{equation}
	For $p=2$, problem (\ref{q}) reduces to the fractional Laplacian problem. we refer to the series of papers of Servadei and Valdinoci \cite{Sevadei2012,Servadei2015,Servadei2013}.
	
	When the nonlinearity is discontinuous, the second author study problem (\ref{q}) and prove the existence and multiplicity results for the following problem
	\begin{equation*}
		\left\{\begin{array}{ll} (-\Delta )^s u = f(u)H(u-\mu) & \quad \mbox{in }\ \Omega\\[0.3cm]
			u= 0        & \quad \mbox{on }\ \mathbb{R}^N \setminus \Omega,
		\end{array}
		\right.
	\end{equation*}
	where $H$ is the Heaviside function, $f$ is a given function and $\mu>0$.
	
	For $p\neq2$ and when $f$ is regular nonlinearity, there is a large number of papers treated problem (\ref{q}) using different techniques. we refer the reader to the papers \cite{Ambrosio2017,Ambrosio2018,Ambrosio2018I,Caffarelli2012,Iannizzotto2016,Mosconi2016,Sevadei2012,Xiang2016} and to the references therein.
	
	When $f$ is discontinuous nonlinearity with respect to $u$, problem (\ref{q}) was treated by the authors in \cite{Achour}. More precisely, $f$ is given by
	$$f(x,u) = m(x)\Sum_{i=1}^{n}~H(u-\mu_i),$$
	for some  $\mu_{i}>0$ verifying the condition
	$$ \mu_1<\mu_2<\cdots<\mu_n,~\text{ for } n \geq 1$$
	and $m \in L^{\infty}(\Omega)$ changes sign.
	
	The authors prove the existence and multiplicity result via the nonsmooth critical point theory.
	
	In this paper, we are interested to study the existence of weak solutions for problem (\ref{q}) with singular absorption term of the following types
	 \begin{equation}
		\label{P}
		\left\{\begin{array}{ll} (-\Delta )_p^s u + \dfrac{|u|^{p-2}u}{|x|^{sp}} = \lambda f(x,u) & \quad \mbox{in }\ \Omega\\[0.3cm]
			u= 0        & \quad \mbox{on }\ \mathbb{R}^N \setminus \Omega,
		\end{array}
		\right.
	\end{equation}
	and
	\begin{equation}
		\label{P-}
		\left\{\begin{array}{ll} (-\Delta )_p^s u - \dfrac{|u|^{p-2}u}{|x|^{sp}} = \lambda f(x,u) & \quad \mbox{in }\ \Omega\\[0.3cm]
			u= 0        & \quad \mbox{on }\ \mathbb{R}^N \setminus \Omega,
		\end{array}
		\right.
	\end{equation}
	where $\lambda$ is a positive parameter, $\Omega$ is a bounded domain in $\mathbb{R}^N$, $p>1$, $s\in(0,1)$, $(N>ps)$ containing the origin and with smooth boundary, $f : \Omega \times \mathbb{R} \longrightarrow \mathbb{R}$ is a Carathéodory function satisfying the following subcritical growth condition
	\begin{enumerate}[label=\textbf{(H1)}]
		\item $\left|f(x,t)\right| \leq \alpha + \beta |t|^{q-1}$, ~~~~$\forall (x,t) \in \Omega \times \mathbb{R}$
		\label{H1}
	\end{enumerate}
	for some non-negative constant $\alpha,\beta$ and $q \in ]p,p^{*}[$, where the fractional critical exponent be defined as
	$$p^{*} \coloneqq \left\{\begin{array}{ll} \frac{Np}{N-sp} & \quad \mbox{if }\ N > ps\\[0.1cm]
		\infty          & \quad \mbox{if }\ N \leq ps.
	\end{array}
	\right.$$
	In the local setting $(s=1)$, M. Khodabakhshi and al. \cite{Khodabakhshi2016} studied the existence of solutions to the problem (\ref{P}) which was motivated by the work of Ferrara and Bisci \cite{Ferrera2014}. They studied the existence of at least one nontrivial solution of the following elliptic problem
	\begin{equation}
		\label{F}
		\left\{\begin{array}{ll} -\Delta_p u = \mu \dfrac{|u|^{p-2}u}{|x|^{p}}+\lambda f(x,u) & \quad \mbox{in }\ \Omega\\[0.3cm]
			u= 0        & \quad \mbox{on }\ \partial \Omega,
		\end{array}
		\right.
	\end{equation}
	where $-\Delta_p u \coloneqq \text{div } (\vert \nabla u\vert^{p-2}\nabla u)$ denotes the $p$-Laplacian, $\lambda,\mu>0$.
	
	Note also the work of Khodabakhshi and Hadjian in \cite{Khodabakhshi2018} when the authors prove the existences of three weak solutions of the following problem
	\begin{equation}
		\label{P2}
		\left\{\begin{array}{ll} -\Delta_p u +\dfrac{|u|^{p-2}u}{|x|^{p}} = \lambda f(x,u) + \mu g(x,u) & \quad \mbox{in }\ \Omega\\[0.3cm]
			u= 0        & \quad \mbox{on }\ \partial \Omega,
		\end{array}
		\right.
	\end{equation}
	where $f$ and $g$ are Carathéodory functions.
	
	Hence, in this paper, we consider the nonlocal counter part of problem (\ref{P2}) when $\mu=0$ and the counter part of problem (\ref{F}). Using a variational structure of our problems and based on a version of critical point theorem contained in \cite{Bonanno2012,Ricceri2000} (See respectively Theorem \ref{th2} and Theorem \ref{th22}), we prove the existence of two weak solutions of problem (\ref{P}) and we show that the problem (\ref{P-}) admits at least one nontrivial solution .
	
	The paper is organized as follows. In Section \ref{sec2}, we recall some basic definitions. Section \ref{sec3}  is denoted to state and prove our main results and finally some useful comments related to problems (\ref{P}) and (\ref{P-}) are given.
	
	\section{Preliminaries.}
	\label{sec2}
	In this section, we provide the variational setting for the problem (\ref{P}), jointly with some preliminary results for the fractional $p$-Laplacian. Also, we recall certain definitions and essential results on the critical point theory.
	\subsection{Variational Formulation of the Problem.}
	Let $\Omega \subset \mathbb{R}^N$ be a bounded domain containing the origin and with smooth boundary $\partial \Omega$, for $p\in(1,\infty)$ and $s\in(0,1)$, we define the fractional Sobolev space
	$$W^{s,p}(\mathbb{R}^N)\coloneqq \left\lbrace  u \in L^{p}(\mathbb{R}^N)~:~\dfrac{\vert u(x)-u(y)\vert}{\vert x-y \vert^{\frac{N}{p}+s}} \in L^{p}(\mathbb{R}^{2N})\right\rbrace ,$$
	endowed with the norm
	$$\Vert u\Vert_{W^{s,p}(\mathbb{R}^N)} \coloneqq \left(  \Vert u\Vert_{L^{p}(\mathbb{R}^N)}^{p} + \left[ u\right]_{W^{s,p}(\mathbb{R}^N)}^{p}\right)^{\frac{1}{p}} ,$$
	where the following term is the Gagliardo semi-norm of u defined as
	$$\left[ u\right]_{W^{s,p}(\mathbb{R}^N)} \coloneqq \left( \int_{\mathbb{R}^{2N}}\dfrac{\vert u(x)-u(y) \vert^{p}}{\vert x-y \vert^{N+sp}}~dxdy\right)^{\frac{1}{p}},$$
	We shall work in the following closed linear subspace
	$$W_{0} \coloneqq \left\lbrace u \in W^{s,p}(\mathbb{R}^N)~:~ u=0 ~\mbox{ a.e. in }~ \mathbb{R}^N \setminus \Omega \right\rbrace ,$$
	which can be equivalently renormed by
	$$\Vert u\Vert_{W_{0}} \coloneqq \left[ u\right]_{W^{s,p}(\mathbb{R}^N)}.$$
	
	Let $1<ps<N$, then there exist a positive constant $c_{H}$ such that we present the fractional Hardy inequality, which says that
	\begin{equation}
		\label{FHardy}
		\int_{\mathbb{R}^{N}}\dfrac{\vert u(x)\vert^{p}}{\vert x \vert^{sp}}~dx \leq c_{H} \int_{\mathbb{R}^{2N}}\dfrac{\vert u(x)-u(y) \vert^{p}}{\vert x-y \vert^{N+sp}}~dxdy, ~~~~\forall u \in W_{0}.
	\end{equation}
	For more details, see \cite{Frank2008}.\\
	Now, let $\mathbb{C}_c^{\infty}(\mathbb{R}^N)$ the space of functions on $\mathbb{R}$ that are infinitely differentiable and have compact support contained in $\mathbb{R}^{N}$. If we denote by $\mathbb{D}^{s,p}(\mathbb{R}^N)$ the closure of $\mathbb{C}_c^{\infty}(\mathbb{R}^N)$ with respect to the Gagliardo semi-norm $\left[ u\right]_{W^{s,p}(\mathbb{R}^N)},$ then for $sp<N,$ we define the Sobolev constant $S_p$ by
	\begin{equation}
		\label{sobolev}
		S_p=\Inf_{v\in \mathbb{D}^{s,p}(\mathbb{R}^N)}\frac{\left[ u\right]_{W^{s,p}}^p}{\|v\|_{L^{p_s^*}}^p},
	\end{equation}
	where $p_{s}^{*}=\frac{Np}{N-sp}$ is the critical Sobolev exponent.\\\\
	Also, It is well-known that $(W_{0}, \Vert \cdot\Vert_{W_{0}})$ is a uniformly convex (i.e. reflexive) Banach space, continuously embedded into $L^{p}(\Omega)$ for all $p \in [1,p^{*}]$ and compactly injected in $L^{p}(\Omega)$ for all $p \in [1,p^{*})$, there exist a $c_p > 0$ which is the best constant of the embedding, such that
	\begin{equation}
		\label{lp}
		||u||_{L^{p}(\Omega)} \leq c_p ||u||_{W_{0}}, ~~~~\forall u \in W_{0}.
	\end{equation}
	Next, let us denote by $(W^{*}_{0}, \Vert \cdot\Vert_{W^{*}_{0}})$ the dual space of $(W_{0}, \Vert \cdot\Vert_{W_{0}})$ and we define the nonlinear operators $A_{p,m} : W_{0} \rightarrow W^{*}_{0}$ as
	$$\langle A_{p,m} (u), v\rangle = \displaystyle \int_{\mathbb{R}^{2N}}\frac{\vert u(x)-u(y) \vert^{p-2}(u(x)-u(y)) (v(x)-v(y))}{\vert x-y \vert^{N+sp}}~dxdy +m \int_{\Omega}\frac{\vert u(x)\vert^{p-2}}{\vert x \vert^{sp}} u(x)v(x)~dx,$$
	where $\langle \cdot, \cdot \rangle$ denotes the inner product on $W_{0}$ and $m=1$ or $m=-1.$
	
	\begin{lem}
		\label{lem1}
		For $u,v \in W_{0}$ and there exist a constant $k \geq 1$, then the nonlinear operators $A_{p,m}$ are well defined and verify the following:
		$$\langle A_{p,m} (u), v\rangle \leq k \Vert u\Vert_{W_{0}}^{p-1} \Vert v\Vert_{W_{0}}$$
		and
		$$\Vert A_{p,m} (u)\Vert_{W^{*}_{0}} \leq k \Vert u\Vert_{W_{0}}^{p-1}.$$
	\end{lem}
	
	\pr For all $u,v \in W_{0}$, we have
	\begin{align*}
		\langle A_{p,m} (u), v\rangle &= \displaystyle \int_{\mathbb{R}^{2N}}\dfrac{\vert u(x)-u(y) \vert^{p-2}(u(x)-u(y)) (v(x)-v(y))}{\vert x-y \vert^{N+sp}}~dxdy +m \int_{\Omega}\dfrac{\vert u(x)\vert^{p-2}}{\vert x \vert^{sp}} u(x)v(x)~dx\\
		&\leq \displaystyle \int_{\mathbb{R}^{2N}}\dfrac{\vert u(x)-u(y) \vert^{p-2}|u(x)-u(y)| |v(x)-v(y)|}{\vert x-y \vert^{(N+sp)\left( \frac{p-1+1}{p}\right)}}~dxdy +m \int_{\Omega}\dfrac{\vert u(x)\vert^{p-2}}{\vert x \vert^{(sp)\left( \frac{p-1+1}{p}\right)}} |u(x)||v(x)|~dx\\
		&\leq \displaystyle \int_{\mathbb{R}^{2N}}\dfrac{\vert u(x)-u(y) \vert^{p-1}}{\vert x-y \vert^{(N+sp)\left( \frac{p-1}{p}\right)}} \dfrac{|v(x)-v(y)|}{\vert x-y \vert^{\left( \frac{N+sp}{p}\right)}}~dxdy +m \int_{\Omega}\dfrac{\vert u(x)\vert^{p-1}}{\vert x \vert^{(sp)\left( \frac{p-1}{p}\right)}} \dfrac{|v(x)|}{\vert x \vert^{\frac{ps}{p}}}~dx.
	\end{align*}
	Then, by Holder inequality we get
	\begin{align*}
		\langle A_{p,m} (u), v\rangle &\leq \left( \displaystyle \int_{\mathbb{R}^{2N}}\dfrac{\vert u(x)-u(y) \vert^{p}}{\vert x-y \vert^{N+sp}}~dxdy \right)^{\frac{p-1}{p}} \left( \int_{\mathbb{R}^{2N}} \dfrac{|v(x)-v(y)|}{\vert x-y \vert^{N+sp}}~dxdy\right)^{\frac{1}{p}}\\
		&~~~~~~+m \left( \int_{\Omega}\dfrac{\vert u(x)\vert^{p}}{\vert x \vert^{sp}}~dx \right)^{\frac{p-1}{p}} \left( \int_{\Omega}\dfrac{\vert v(x)\vert^{p}}{\vert x \vert^{sp}}~dx \right)^{\frac{1}{p}}.
	\end{align*}
	For any $\beta \in (0,1)$ and $a,b,c,d > 0$, we use the following inequality
	$$a^{\beta}c^{1-\beta} +m ~b^{\beta} d^{1-\beta} \leq (a+m ~b)^{\beta}(c+m ~d)^{1-\beta},$$
	by setting $\beta = \frac{p-1}{p}$ and
	\begin{equation*}
		\begin{array}{l}
			a = \displaystyle \int_{\mathbb{R}^{2N}}\dfrac{\vert u(x)-u(y) \vert^{p}}{\vert x-y \vert^{N+sp}}~dxdy\\[0.4cm]
			b = \displaystyle \int_{\Omega}\dfrac{\vert u(x)\vert^{p}}{\vert x \vert^{sp}}~dx\\[0.4cm]
			c = \displaystyle \int_{\mathbb{R}^{2N}} \dfrac{|v(x)-v(y)|}{\vert x-y \vert^{N+sp}}~dxdy\\[0.4cm]
			d = \displaystyle \int_{\Omega}\dfrac{\vert v(x)\vert^{p}}{\vert x \vert^{sp}}~dx,
		\end{array}
	\end{equation*}
	we can deduce that
	\begin{align*}
		\langle A_{p,m} (u), v\rangle &\leq \left( \displaystyle \int_{\mathbb{R}^{2N}}\dfrac{\vert u(x)-u(y) \vert^{p}}{\vert x-y \vert^{N+sp}}~dxdy +m \int_{\Omega}\dfrac{\vert u(x)\vert^{p}}{\vert x \vert^{sp}}~dx \right)^{\frac{p-1}{p}}\\ &~~~~~~ \times \left( \int_{\mathbb{R}^{2N}} \dfrac{|v(x)-v(y)|}{\vert x-y \vert^{N+sp}}~dxdy +m \int_{\Omega}\dfrac{\vert v(x)\vert^{p}}{\vert x \vert^{sp}}~dx \right)^{\frac{1}{p}}.
	\end{align*}
	Then, according to the fractional Hardy inequality (\ref{FHardy}) we get
	\begin{align*}
		\langle A_{p,1} (u), v\rangle &\leq \left( \displaystyle \int_{\mathbb{R}^{2N}}\dfrac{\vert u(x)-u(y) \vert^{p}}{\vert x-y \vert^{N+sp}}~dxdy + c_{H} \int_{\mathbb{R}^{2N}}\dfrac{\vert u(x)-u(y) \vert^{p}}{\vert x-y \vert^{N+sp}}~dxdy \right)^{\frac{p-1}{p}}\\ &~~~~~~ \times \left( \int_{\mathbb{R}^{2N}} \dfrac{|v(x)-v(y)|}{\vert x-y \vert^{N+sp}}~dxdy + c_{H} \int_{\mathbb{R}^{2N}} \dfrac{|v(x)-v(y)|}{\vert x-y \vert^{N+sp}}~dxdy \right)^{\frac{1}{p}}\\
		&\leq (c_{H}+1) \Vert u\Vert_{W_{0}}^{p-1} \Vert v\Vert_{W_{0}}\\
		&\leq k \Vert u\Vert_{W_{0}}^{p-1} \Vert v\Vert_{W_{0}}\\
		&< + \infty.
	\end{align*}
	Where $k=c_{H}+1$ and $k \geq 1$, we have
	\begin{align*}
		\langle A_{p,-1} (u), v\rangle &\leq \left( \displaystyle \int_{\mathbb{R}^{2N}}\dfrac{\vert u(x)-u(y) \vert^{p}}{\vert x-y \vert^{N+sp}}~dxdy - \int_{\Omega}\dfrac{\vert u(x)\vert^{p}}{\vert x \vert^{sp}}~dx \right)^{\frac{p-1}{p}}\\ &~~~~~~ \times \left( \int_{\mathbb{R}^{2N}} \dfrac{|v(x)-v(y)|}{\vert x-y \vert^{N+sp}}~dxdy - \int_{\Omega}\dfrac{\vert v(x)\vert^{p}}{\vert x \vert^{sp}}~dx \right)^{\frac{1}{p}}\\
		&\leq  \Vert u\Vert_{W_{0}}^{p-1} \Vert v\Vert_{W_{0}}\\
		&< + \infty.
	\end{align*}
	Moreover, respectively we have
	$$\Vert A_{p,1} (u)\Vert_{W^{*}_{0}} \leq k \Vert u\Vert_{W_{0}}^{p-1}$$
	and
	$$\Vert A_{p,-1} (u)\Vert_{W^{*}_{0}} \leq \Vert u\Vert_{W_{0}}^{p-1}.$$ \hspace*{\fill} $\square$
	
	\begin{lem}
		\label{lem2}
		For all $u,v \in W_{0}$ and there exist a constant $\mu \leq 1$, the operators $A_{p,m}$ satisfy the following inequalities :
		$$\langle A_{p,1} (u)- A_{p,1} (v), u-v\rangle \geq k \left( \Vert u\Vert_{W_{0}}^{p-1} - \Vert v\Vert_{W_{0}}^{p-1}\right) \left(\Vert u\Vert_{W_{0}} - \Vert v\Vert_{W_{0}} \right)$$
		and
		$$\langle A_{p,-1} (u)- A_{p,-1} (v), u-v\rangle \geq \mu \left( \Vert u\Vert_{W_{0}}^{p-1} - \Vert v\Vert_{W_{0}}^{p-1}\right) \left(\Vert u\Vert_{W_{0}} - \Vert v\Vert_{W_{0}} \right).$$
	\end{lem}
	\pr By direct computation, we have \\
	\resizebox{1\linewidth}{!}{\begin{minipage}{\linewidth}
			\begin{align*}
				\langle A_{p,m} (u)- A_{p,m} (v), u-v\rangle &= \langle A_{p,m} (u), u-v\rangle - \langle A_{p,m} (v), u-v\rangle \\
				&= \displaystyle \int_{\mathbb{R}^{2N}}\dfrac{\vert u(x)-u(y) \vert^{p-2}(u(x)-u(y)) ((u-v)(x)-(u-v) (y))}{\vert x-y \vert^{N+sp}}~dxdy\\
				&~~~~~+m \int_{\Omega}\dfrac{\vert u(x)\vert^{p-2}}{\vert x \vert^{sp}} u(x)(u-v)(x)~dx\\
				&~~~~~-\displaystyle \int_{\mathbb{R}^{2N}}\dfrac{\vert v(x)-v(y) \vert^{p-2}(v(x)-v(y)) ((u-v)(x)-(u-v) (y))}{\vert x-y \vert^{N+sp}}~dxdy\\
				&~~~~~+m \int_{\Omega}\dfrac{\vert v(x)\vert^{p-2}}{\vert x \vert^{sp}} v(x)(u-v)(x)~dx\\
				&= \displaystyle \int_{\mathbb{R}^{2N}}\dfrac{\vert u(x)-u(y) \vert^{p}}{\vert x-y \vert^{N+sp}}~dxdy +m \int_{\Omega}\dfrac{\vert u\vert^{p}}{\vert x \vert^{sp}}~dx + \displaystyle \int_{\mathbb{R}^{2N}}\dfrac{\vert v(x)-v(y) \vert^{p}}{\vert x-y \vert^{N+sp}}~dxdy +m \int_{\Omega}\dfrac{\vert v\vert^{p}}{\vert x \vert^{sp}}~dx \\
				&- \displaystyle \int_{\mathbb{R}^{2N}}\dfrac{\vert u(x)-u(y) \vert^{p-2}(u(x)-u(y)) (v(x)-v(y))}{\vert x-y \vert^{N+sp}}~dxdy +m \int_{\Omega}\dfrac{\vert u(x)\vert^{p-2}}{\vert x \vert^{sp}} u(x)v(x)~dx\\
				&- \displaystyle \int_{\mathbb{R}^{2N}}\dfrac{\vert v(x)-v(y) \vert^{p-2}(v(x)-v(y)) (u(x)-u(y))}{\vert x-y \vert^{N+sp}}~dxdy +m \int_{\Omega}\dfrac{\vert v(x)\vert^{p-2}}{\vert x \vert^{sp}} v(x)u(x)~dx\\
				&= \displaystyle \int_{\mathbb{R}^{2N}}\dfrac{\vert u(x)-u(y) \vert^{p}}{\vert x-y \vert^{N+sp}}~dxdy +m \int_{\Omega}\dfrac{\vert u\vert^{p}}{\vert x \vert^{sp}}~dx + \displaystyle \int_{\mathbb{R}^{2N}}\dfrac{\vert v(x)-v(y) \vert^{p}}{\vert x-y \vert^{N+sp}}~dxdy +m \int_{\Omega}\dfrac{\vert v\vert^{p}}{\vert x \vert^{sp}}~dx\\
				&~~~~~- \langle A_{p,m} (u), v\rangle - \langle A_{p,m} (v), u\rangle.\\
			\end{align*}
	\end{minipage}}\\
	Furthermore, by Lemma \ref{lem1}, there exist $k \geq 1$ such that
	\begin{align*}
		\langle A_{p,1} (u)- A_{p,1} (v), u-v\rangle &\geq \Vert u\Vert_{W_{0}}^{p} + \Vert v\Vert_{W_{0}}^{p}  - \langle A_{p,1} (u), v\rangle - \langle A_{p,1} (v), u\rangle \\
		&\geq \Vert u\Vert_{W_{0}}^{p} + \Vert v\Vert_{W_{0}}^{p} - k \Vert u\Vert_{W_{0}}^{p-1} \Vert v\Vert_{W_{0}} - k \Vert v\Vert_{W_{0}}^{p-1} \Vert u\Vert_{W_{0}}\\
		&= - k \left( - \frac{1}{k}\Vert u\Vert_{W_{0}}^{p} - \frac{1}{k} \Vert v\Vert_{W_{0}}^{p} + \Vert u\Vert_{W_{0}}^{p-1} \Vert v\Vert_{W_{0}} + \Vert v\Vert_{W_{0}}^{p-1} \Vert u\Vert_{W_{0}} \right)\\
		&\geq - k \left( -\Vert u\Vert_{W_{0}}^{p} - \Vert v\Vert_{W_{0}}^{p} + \Vert u\Vert_{W_{0}}^{p-1} \Vert v\Vert_{W_{0}} + \Vert v\Vert_{W_{0}}^{p-1} \Vert u\Vert_{W_{0}} \right)\\
		&\geq k \left( \Vert u\Vert_{W_{0}}^{p-1} - \Vert v\Vert_{W_{0}}^{p-1}\right)
		\left(\Vert u\Vert_{W_{0}} - \Vert v\Vert_{W_{0}} \right).
	\end{align*}
	Likewise, according to the fractional Hardy inequality (\ref{FHardy}) and Lemma \ref{lem1} we get
	\begin{align*}
		\langle A_{p,-1} (u)- A_{p,-1} (v), u-v\rangle &\geq (1-c_{H})\Vert u\Vert_{W_{0}}^{p} + (1-c_{H})\Vert v\Vert_{W_{0}}^{p}  - \langle A_{p,-1} (u), v\rangle - \langle A_{p,-1} (v), u\rangle\\
		&\geq \mu \Vert u\Vert_{W_{0}}^{p} + \mu \Vert v\Vert_{W_{0}}^{p} - \Vert u\Vert_{W_{0}}^{p-1} \Vert v\Vert_{W_{0}} - \Vert v\Vert_{W_{0}}^{p-1} \Vert u\Vert_{W_{0}}\\
		&= -\mu \left( -\Vert u\Vert_{W_{0}}^{p} - \Vert v\Vert_{W_{0}}^{p} + \dfrac{1}{\mu} \Vert u\Vert_{W_{0}}^{p-1} \Vert v\Vert_{W_{0}} + \dfrac{1}{\mu} \Vert v\Vert_{W_{0}}^{p-1} \Vert u\Vert_{W_{0}} \right)\\
		&\geq - \mu \left( -\Vert u\Vert_{W_{0}}^{p} - \Vert v\Vert_{W_{0}}^{p} + \Vert u\Vert_{W_{0}}^{p-1} \Vert v\Vert_{W_{0}} + \Vert v\Vert_{W_{0}}^{p-1} \Vert u\Vert_{W_{0}} \right)\\
		&\geq \mu \left( \Vert u\Vert_{W_{0}}^{p-1} - \Vert v\Vert_{W_{0}}^{p-1}\right)
		\left(\Vert u\Vert_{W_{0}} - \Vert v\Vert_{W_{0}} \right),
	\end{align*}
	where $\mu=1-c_{H}$ and $\mu\leq1$.\hspace*{\fill} $\square$\\
	
	\begin{prop} ~\smallskip The nonlinear operators $A_{p,m}$ have the following properties :
		\label{lemp}
		\begin{enumerate}[label= \arabic*)]
			\item $A_{p,m} : W_{0} \rightarrow W^{*}_{0}$ are a continuous, bounded and strictly monotone operators, i.e: if $\langle A_{p,m} (u)-A_{p,m} (v), u-v\rangle > 0$, $\forall u \neq v$. \label{(1)}
			
			\item $A_{p,m}$ are a mappings of type $(S)$, i.e: if $u_{n} \rightarrow u$  weakly in $W_{0}$ and ~$\displaystyle\limsup_{n \rightarrow \infty}\langle A_{p,m} (u_{n})-A_{p,m} (u), u_{n}-u\rangle \leq 0$, then $u_{n} \rightarrow u$ strongly in  $W_{0}$. \label{(2)}
			
			\item $A_{p,m} : W_{0} \rightarrow W^{*}_{0}$ are homomorphisms.
		\end{enumerate}
	\end{prop}
	\pr  ~\smallskip
	
	\noindent 1) By Lemma \ref{lem1}, there exist a constant $k \geq 1$ such that
	$$\left| \langle A_{p,m} (u), v\rangle\right| \leq k \Vert u\Vert_{W_{0}}^{p-1} \Vert v\Vert_{W_{0}},~\forall u,v \in W_{0}$$
	From this inequality, it is obvious that $A_{p,m}$ are continuous and bounded.\\
	
	Now, using the well-known Simon inequality (see \cite[formula (2.2)]{Simon1978}): for all $\xi,\eta \in \mathbb{R}^{N}$, for $p>1$ there exists a positive constant $C_{p}$ such that
	\begin{equation}
		\label{sim_ineq}
		C_{p}\langle |\xi|^{p-2}\xi-|\eta|^{p-2}\eta, \xi-\eta\rangle \geq \left\{\begin{array}{ll}
			\left| \xi-\eta\right|^{p} & \quad \mbox{if }\ p\geq2,\\
			\left| \xi-\eta\right|^{2}\left( |\xi|^{p}+|\eta|^{p}\right)^{(p-2)/p} & \quad \mbox{if }\ 1<p<2.
		\end{array}
		\right.
	\end{equation}
	Then, using Lemma \ref{lem2} and applying inequality (\ref{sim_ineq}) for any $u,v \in W_{0}$, with $u \neq v$, we have if $p\geq2$
	$$C_{p}\langle A_{p,1} (u)-A_{p,1} (v), u-v\rangle \geq k \left( \Vert u\Vert_{W_{0}}^{p-1} - \Vert v\Vert_{W_{0}}^{p-1}\right) \left(\Vert u\Vert_{W_{0}} - \Vert v\Vert_{W_{0}} \right) \geq 0$$
	equally,
	$$C_{p}\langle A_{p,-1} (u)-A_{p,-1} (v), u-v\rangle \geq \mu \left( \Vert u\Vert_{W_{0}}^{p-1} - \Vert v\Vert_{W_{0}}^{p-1}\right) \left(\Vert u\Vert_{W_{0}} - \Vert v\Vert_{W_{0}} \right) \geq 0.$$
	
	and if $1<p<2$, we have
	$$C_{p}^{p/2} \left[ \langle A_{p,1} (u)-A_{p,1} (v), u-v\rangle\right]^{p/2}\left( \Vert u\Vert_{W_{0}}^{p} - \Vert v\Vert_{W_{0}}^{p}\right)^{(2-p)/2} \geq k \left( \Vert u\Vert_{W_{0}}^{p-1} - \Vert v\Vert_{W_{0}}^{p-1}\right) \left(\Vert u\Vert_{W_{0}} - \Vert v\Vert_{W_{0}} \right),$$
	likewise,
	$$C_{p}^{p/2} \left[ \langle A_{p,-1} (u)-A_{p,-1} (v), u-v\rangle\right]^{p/2}\left( \Vert u\Vert_{W_{0}}^{p} - \Vert v\Vert_{W_{0}}^{p}\right)^{(2-p)/2} \geq \mu \left( \Vert u\Vert_{W_{0}}^{p-1} - \Vert v\Vert_{W_{0}}^{p-1}\right) \left(\Vert u\Vert_{W_{0}} - \Vert v\Vert_{W_{0}} \right),$$
	thus,
	$$C\left[ \langle A_{p,1} (u)-A_{p,-1} (v), u-v\rangle\right]^{p/2} \geq k \left( \Vert u\Vert_{W_{0}}^{p-1} - \Vert v\Vert_{W_{0}}^{p-1}\right) \left(\Vert u\Vert_{W_{0}} - \Vert v\Vert_{W_{0}} \right)\geq 0$$
	and
	$$C\left[ \langle A_{p,-1} (u)-A_{p,-1} (v), u-v\rangle\right]^{p/2} \geq \mu \left( \Vert u\Vert_{W_{0}}^{p-1} - \Vert v\Vert_{W_{0}}^{p-1}\right) \left(\Vert u\Vert_{W_{0}} - \Vert v\Vert_{W_{0}} \right)\geq 0,$$
	where $C>0$ is a constant, which lead us to conclude that $A_{p,m}$ are strictly monotone.\\
	
	\noindent 2) Since $W_{0}$ is a reflexive Banach space, it is isometrically isomorphic to a locally uniformly convex space. So as it was already proved, weak convergence and norm convergence imply strong convergence. Therefore we only need to show that $\Vert u_{n}\Vert_{W_{0}} \rightarrow \Vert u\Vert_{W_{0}}$.\\
	
	\noindent further, we have that if $u_{n} \rightarrow u$  weakly in $W_{0}$ and $$\displaystyle\limsup_{n \rightarrow \infty}\langle A_{p,m} (u_{n})-A_{p,m} (u), u_{n}-u\rangle \leq 0.$$
	\noindent Then
	$$\displaystyle \lim_{n \rightarrow +\infty} \langle A_{p,m} (u_{n})-A_{p,m} (u), u_{n}-u\rangle = \lim_{n \rightarrow +\infty} \langle A_{p,m} (u_{n}), u_{n}-u\rangle - \langle A_{p,m} (u), u_{n}-u\rangle = 0.$$
	By Lemma \ref{lem2} and \ref{(1)}, we have
	$$\langle A_{p,1} (u_{n})-A_{p,1} (u), u_{n}-u\rangle \geq k \left( \Vert u_{n}\Vert_{W_{0}}^{p-1} - \Vert u\Vert_{W_{0}}^{p-1}\right) \left(\Vert u_{n}\Vert_{W_{0}} - \Vert u\Vert_{W_{0}} \right) \geq 0$$
	and
	$$\langle A_{p,-1} (u_{n})-A_{p,-1} (u), u_{n}-u\rangle \geq \mu \left( \Vert u_{n}\Vert_{W_{0}}^{p-1} - \Vert u\Vert_{W_{0}}^{p-1}\right) \left(\Vert u_{n}\Vert_{W_{0}} - \Vert u\Vert_{W_{0}} \right) \geq 0.$$
	Which imply that $\Vert u_{n}\Vert_{W_{0}} \rightarrow \Vert u\Vert_{W_{0}}$ as $n \rightarrow + \infty$ and lead us to conclude that $u_{n} \rightarrow u$ strongly in  $W_{0}$ as $n \rightarrow + \infty$.\\
	
	\noindent 3) By \ref{(1)}, we know that $A_{p,m}$ are strictly monotone, which implicate that $A_{p,m}$ are injective.
	Also, according to Lemma \ref{lem1}, we have
	\begin{align*}
		\lim_{\Vert u\Vert_{W_{0}} \rightarrow +\infty} \dfrac{\langle A_{p,m} (u), u\rangle}{\Vert u\Vert_{W_{0}}} &= \lim_{\Vert u\Vert_{W_{0}} \rightarrow +\infty} \dfrac{\Vert u\Vert_{W_{0}}^{p} +m \left\Vert \frac{u}{|x|^{s}}\right\Vert _{L^{p}(\Omega)}^{p} }{\Vert u\Vert_{W_{0}}}\\
		&= \lim_{\Vert u\Vert_{W_{0}} \rightarrow +\infty} \Vert u\Vert_{W_{0}}^{p-1} +m  \left\Vert \frac{u}{|x|^{s}}\right\Vert _{L^{p}(\Omega)}^{p} \Vert u\Vert_{W_{0}}^{-1}\\
		&= \lim_{\Vert u\Vert_{W_{0}} \rightarrow +\infty} \Vert u\Vert_{W_{0}}^{p-1}\\
		&= + \infty.
	\end{align*}
	Thanks to $1<p<\frac{N}{s}$, hence $A_{p,m}$ are Coercive on $W_{0}$. Since $A_{p,m}$ are continuous and bounded by \ref{(1)}. Then, by the Minty-Browder Theorem (see \cite[Theorem 26.A]{Zeidler1990B}), we conclude that $A_{p,m}$ are a surjections.\\
	Thus, $A_{p,m}$ have inverse mappings $A_{p,m}^{-1} : W_{0}^{*} \rightarrow W_{0}$. Therefore, the continuity of $A_{p,m}^{-1}$ is sufficient to ensure $A_{p,m}$ to be homeomorphisms.\\
	Assume, $g_{n},g \in W_{0}$ with $g_{n} \rightarrow g$ in $W_{0}$. Let $u_{n}=A_{p,m}^{-1}(g_{n})$ and $u=A_{p,m}^{-1}(g)$. Then $A_{p,m}(u_{n})=g_{n}$ and $A_{p,m}(u)=g$. Clearly, $\{ u_{n}\}$ is bounded in $W_{0}$. Thus there exist $u_{0} \in W_{0}$ and a subsequence of $\{ u_{n}\}$ still denoted by $\{ u_{n}\}$ such that $u_{n} \rightharpoonup u_{0}$ since $g_{n} \rightarrow g$, we have
	$$\lim_{n \rightarrow +\infty} \langle A_{p,m} (u_{n})-A_{p,m} (u_{0}), u_{n}-u_{0}\rangle = \lim_{n \rightarrow +\infty} \langle g_{n}, u_{n}-u_{0}\rangle =0.$$
	In view of $A_{p,m}$ satisfying the $(S)$ condition by \ref{(2)}, we get $u_{n} \rightarrow u_{0}$ in $W_{0}$. Moreover, $u=u_{0}$ a.e. in $\Omega$. Hence, $u_{n} \rightarrow u$ in $W_{0}$, so that $A_{p,m}^{-1}$ are continuous. \hspace*{\fill} $\square$\\
	
	Initially, let us introduce the energy functionals $E_{\lambda}^{p,m} : W_{0} \rightarrow \mathbb{R}$ associated with problem (\ref{P}) as
	$$E_{\lambda}^{p,m}(u) \coloneqq \Phi_{p,m} (u) - \lambda \Psi (u), ~~~~\forall u \in W_{0},$$
	where
	$$\Phi_{p,m} (u) \coloneqq \dfrac{1}{p} \left( \displaystyle \int_{\mathbb{R}^{2N}}\dfrac{\vert u(x)-u(y) \vert^{p}}{\vert x-y \vert^{N+sp}}~dxdy +m \int_{\Omega}\dfrac{\vert u(x)\vert^{p}}{\vert x \vert^{sp}}~dx\right)$$
	and
	$$\Psi (u) \coloneqq \int_{\Omega}F(x,u(x))~dx,$$
	where $F(x,t) = \int_{0}^{t}f(x,s)~ds$, for every $(x,t) \in \Omega \times \mathbb{R}$.\\\\
	Moreover, from the fractional Hardy's inequality (\ref{FHardy}), we have
	\begin{align*}
		\Phi_{p,1}(u) &= \dfrac{1}{p} \displaystyle \int_{\mathbb{R}^{2N}}\dfrac{\vert u(x)-u(y) \vert^{p}}{\vert x-y \vert^{N+sp}}~dxdy + \dfrac{1}{p} \int_{\Omega}\dfrac{\vert u(x)\vert^{p}}{\vert x \vert^{sp}}~dx\\
		&\leq  \dfrac{1}{p} \Vert u\Vert^{p}_{W_{0}} + \dfrac{c_{H}}{p} \Vert u\Vert^{p}_{W_{0}}\\
		&\leq \left( \dfrac{c_{H}+1}{p}\right) \Vert u\Vert^{p}_{W_{0}}
	\end{align*}
	and
	\begin{align}
		\label{Heqm}
		\Phi_{p,-1}(u) &= \dfrac{1}{p} \displaystyle \int_{\mathbb{R}^{2N}}\dfrac{\vert u(x)-u(y) \vert^{p}}{\vert x-y \vert^{N+sp}}~dxdy - \dfrac{1}{p} \int_{\Omega}\dfrac{\vert u(x)\vert^{p}}{\vert x \vert^{sp}}~dx\nonumber\\
		&\geq  \dfrac{1}{p} \Vert u\Vert^{p}_{W_{0}} - \dfrac{c_{H}}{p} \int_{\Omega}\dfrac{\vert u(x)\vert^{p}}{\vert x \vert^{sp}}~dx \nonumber\\
		&\geq  \left( \dfrac{1-c_{H}}{p}\right) \Vert u\Vert^{p}_{W_{0}}.
	\end{align}
	Also,
	\begin{equation}
		\label{Heqp}
		\Phi_{p,1}(u) = \dfrac{1}{p} \Vert u\Vert^{p}_{W_{0}} + \dfrac{1}{p} \int_{\Omega}\dfrac{\vert u(x)\vert^{p}}{\vert x \vert^{sp}}~dx \geq  \dfrac{1}{p} \Vert u\Vert^{p}_{W_{0}}
	\end{equation}
	and
	\begin{align*}
		\Phi_{p,-1}(u) &= \dfrac{1}{p} \displaystyle \int_{\mathbb{R}^{2N}}\dfrac{\vert u(x)-u(y) \vert^{p}}{\vert x-y \vert^{N+sp}}~dxdy - \dfrac{1}{p} \int_{\Omega}\dfrac{\vert u(x)\vert^{p}}{\vert x \vert^{sp}}~dx\\
		&\leq  \dfrac{1}{p} \Vert u\Vert^{p}_{W_{0}}.
	\end{align*}
	It follows that
	$$\dfrac{1}{p} \Vert u\Vert^{p}_{W_{0}} \leq \Phi_{p,1}(u) \leq \left( \dfrac{c_{H}+1}{p}\right) \Vert u\Vert^{p}_{W_{0}},~~~~\forall u \in W_{0}$$
	equally,
	$$\left( \dfrac{1-c_{H}}{p}\right) \Vert u\Vert^{p}_{W_{0}} \leq \Phi_{p,-1}(u) \leq \dfrac{1}{p} \Vert u\Vert^{p}_{W_{0}},~~~~\forall u \in W_{0}.$$
	So, $\Phi_{p,m}(u)$ are well defined and coercive in $W_{0}$. \\
	
	Lastly, it's evident that our energy functional $E_{\lambda}^{p,m}$ are well defined and of class $C^1$. The derivative of $E_{\lambda}^{p,m}$ is given by
	$$\langle \left(E_{\lambda}^{p,m}\right)^{'} (u), v\rangle = \langle \Phi_{p,m}^{'} (u), v\rangle - \lambda \langle \Psi^{'} (u), v\rangle,~~\forall u,v \in W_{0}.$$
	
	\begin{defi}
		Fixing the real parameter $\lambda$, a function $u : \Omega \rightarrow \mathbb{R}$ is said to be a weak solution of the problem (\ref{P}) and (\ref{P-}), if $u \in W_{0}$ and
		$$\langle \Phi_{p,m}^{'} (u), v\rangle =\langle A_{p,m}(u), v\rangle = \lambda \langle \Psi^{'} (u), v\rangle, ~~v \in W_{0}$$
		where
		$$ \displaystyle \int_{\mathbb{R}^{2N}}\dfrac{\vert u(x)-u(y) \vert^{p-2}(u(x)-u(y))}{\vert x-y \vert^{N+sp}} (v(x)-v(y))~dxdy +m \int_{\Omega}\dfrac{\vert u(x)\vert^{p-2}}{\vert x \vert^{sp}} u(x)v(x)~dx$$$$ = \lambda \int_{\Omega}f(x,u(x))v(x)~dx,$$
		for every $v \in W_{0}$. Hence, the critical points of $E_{\lambda}^{p,m}$ are exactly the weak solutions of problems (\ref{P}) and (\ref{P-}).
	\end{defi}
	\begin{defi}
		\label{PS}
		Let $X$ be a real Banach space. A Gâteaux differentiable function $E$ satisfies the Palais-Smale condition (in short (PS)-condition), if any sequence $\{ u_{n}\}_{n\in \mathbb{N}}$ such that
		\begin{enumerate}[label= (\alph*)]
			\item $\left\{ E (u_{n})\right\}$ is bounded, i.e.
			$ E (u_{n}) \leq d \coloneqq \displaystyle \sup_{n\in \mathbb{N}} \left\{ E (u_{n})\right\}$,\label{h1}\\
			\vspace{-0.5cm}
			\item $\Vert E^{'} (u_{n})  \Vert_{X^{*}} \longrightarrow 0$ as $n \rightarrow + \infty$,\label{h2}
		\end{enumerate}
		has a convergent subsequence.
	\end{defi}
	Our main tools are the following critical point theorems.
	\begin{theo} \cite[Theorem 3.2.]{Bonanno2012}\\
		\label{th2}
		Let $X$ be a real Banach space and let $\Phi,\Psi : X \rightarrow  \mathbb{R}$ be two continuously Gâteaux differentiable functionals such that $\Phi$ is bounded from below and $\Phi(0)=\Psi(0)=0$. Fix $r>0$ such that $$\sup_{\{ \Phi(u)<r\}} \Psi(u) < +\infty$$ and assume that, for each $$\lambda \in \left]0,\frac{r}{\sup_{\{ \Phi(u)<r\}} \Psi(u)}\right[,$$ the functional $E_{\lambda} \coloneqq \Phi - \lambda \Psi$ satisfies the (PS)-condition and it is unbounded from below. Then, for each $\lambda \in ]0,\frac{r}{\sup_{\{ \Phi(u)<r\}} \Psi(u)}[$, the functional $E_{\lambda}$ admits two distinct critical points.
	\end{theo}
	
	\begin{theo} \cite[Theorem 2.5.]{Ricceri2000}\\
		\label{th22}
		Let $X$ be a real Banach space and let $\Phi,\Psi : X \rightarrow  \mathbb{R}$ be two continuously Gâteaux differentiable functionals such that $\Phi$ is strongly continuous, sequentially weakly lower semicontinious and coercive. Further, assume that $\Psi$ is sequentially weakly upper semicontinious. For every $r>\inf_{X}\Phi$, put
		$$\varphi(r) \coloneqq \displaystyle \inf_{u \in \Phi^{-1}(]-\infty,r[)} \dfrac{\left( \sup_{\{ v \in \Phi^{-1}( ]-\infty,r[ )\}} \Psi(v) \right)-\Psi(u)}{r-\Phi(u)}.$$
		Then, for every $r>\inf_{X}\Phi$ and every $\lambda \in \left]0,1/\varphi(r)\right[$, the restriction $E_{\lambda} \coloneqq \Phi - \lambda \Psi$ to \\$\Phi^{-1}( ]-\infty,r[ )$ admits a global minimum, which is a critical point (local minima) of $E_{\lambda}$ in $X$.
	\end{theo}
    Now, we state the following Theorem which provides the Concentration-Compactness principle in fractional Sobolev spaces. This result is useful for the proof of the semicontinuity property of functionals.
   \begin{theo}
 	\label{thccp}
 	Let $\{u_{n}\}_{n\in \mathbb{N}} \subset W_{0}$ be a weakly convergent sequence with weak limit $u$. Then, there exist two bounded measures $\mu$ and $\nu$, an at
 	most enumerable set of indices $I$ of distinct points $\{x_{i}\}_{i\in I} \subset \mathbb{R}^{N}$, and positive real numbers $\{ \mu_{i}\}_{i\in I}, ~\{ \nu_{i}\}_{i\in I}$, such that the following convergence hold weakly in the sense of measures,
 	\begin{equation}
 		\displaystyle \int_{\mathbb{R}^{N}}\dfrac{\vert u_{n}(x)-u_{n}(y) \vert^{p}}{\vert x-y \vert^{N+sp}}~dydx \rightharpoonup \mu \geq \displaystyle \int_{\mathbb{R}^{N}}\dfrac{\vert u(x)-u(y) \vert^{p}}{\vert x-y \vert^{N+sp}}~dydx + \sum_{i \in I} \mu_{i} \delta_{x_{i}},
 	\end{equation}
 	\begin{equation}
 		\left|u_{n}(x)\right|^{p}dx \rightharpoonup \nu = \left|u(x)\right|^{p}dx + \sum_{i \in I} \nu_{i} \delta_{x_{i}},
 	\end{equation}
 	\begin{equation}
 		\label{renumu}
 		S_{p} ~ \nu_{i} \leq \mu_{i}, ~\forall i \in I,
 	\end{equation}
 	where  $S_{p}$ is the Sobolev constant given by (\ref{sobolev}) and $\delta_{x_{i}}$ denotes the Dirac mass at $x_{i}$. Moreover, if we define
 	\begin{equation}
 		\mu_{\infty} = \displaystyle \lim_{R \to \infty}\limsup_{n \to \infty} \int_{|x|\geq R}\displaystyle \int_{\mathbb{R}^{N}}\dfrac{\vert u_{n}(x)-u_{n}(y) \vert^{p}}{\vert x-y \vert^{N+sp}}~dydx,
 	\end{equation}
 	\begin{equation}
 		\nu_{\infty} = \displaystyle \lim_{R \to \infty}\limsup_{n \to \infty} \int_{|x| \geq R} \left|u_{n}(x)\right|^{p}dx,
 	\end{equation}
 	then
 	\begin{equation}
 		\label{muinf}
 		\displaystyle \limsup_{n \to \infty} \int_{|x| \geq R}\displaystyle \int_{\mathbb{R}^{N}}\dfrac{\vert u_{n}(x)-u_{n}(y) \vert^{p}}{\vert x-y \vert^{N+sp}}~dydx = \mu(\mathbb{R}^{N}) + \mu_{\infty},
 	\end{equation}
 	\begin{equation}
 		\label{nuinf}
 		\displaystyle \limsup_{n \to \infty} \int_{|x| \geq R} \left|u_{n}(x)\right|^{p}dx = \nu(\mathbb{R}^{N}) + \nu_{\infty},
 	\end{equation}
 	\begin{equation}
 		\label{renumuinf}
 		S_{p} ~\nu_{\infty} \leq \mu_{\infty},
 	\end{equation}
 \end{theo}
This result was recently proved in the following mentioned papers of J. F. Bonder and al. in \cite[Theorem 1.1.]{Bonder2018} and K. Ho \& Y. H. Kim in \cite[Theorem 4.1.]{Ho2021}.
	
	\section{Main results}
	\label{sec3}
	In this section, we establish the main results of this paper. the first result is the following
	\begin{theo}
		\label{thp1}
		Let $f : \Omega \times \mathbb{R} \rightarrow \mathbb{R}$ be a Carathéodory function such that condition \ref{H1} holds. Moreover, assume that
		\begin{enumerate}[label=\textbf{(H2)}]
			\item there exist $\theta > p$ and $M>0$ such that
			$$0 < \theta F(x,t)\leq tf(x,t),$$
			for each $x \in \Omega$ and $|t| \geq M$.
			\label{H2}
		\end{enumerate}
		Then, for each $\lambda \in ] 0 ,\overline{\lambda} [$, the following initial case of problem (\ref{P})
		\leqnomode
		\begin{equation*}
			\label{p+}
			\left\{\begin{array}{ll} (-\Delta )_p^s u + \dfrac{|u|^{p-2}u}{|x|^{sp}} = \lambda f(x,u) & \quad \mbox{in }\ \Omega\\[0.3cm]
				u= 0        & \quad \mbox{on }\ \mathbb{R}^N \setminus \Omega,\tag{P$_{+}$}
			\end{array}
			\right.
		\end{equation*}
		\reqnomode
		admits at least two distinct weak solutions, where
		$$\overline{\lambda} \coloneqq \dfrac{q}{qa_{1}c_{1}p^{\frac{1}{p}} + a_{2}c^{q}_{q}p^{\frac{q}{p}}}.$$
	\end{theo}
	To achieve the proof of this result, we may split into a sequence of propositions the proof of Theorem \ref{thp1}
	
\begin{prop}
	\label{prop1}
	For every $\lambda > 0$ and the assumptions \ref{H1}, \ref{H2} hold. Then, $E_{\lambda}^{p,1} = \Phi_{p,1} - \lambda \Psi$ satisfies the Palais-Smale condition.
\end{prop}
\pr To prove that $E_{\lambda}^{p,1}$ satisfies the Palais-Smale condition for every $\lambda > 0$. Namely, we need to show that any sequence of Palais-Smale is bounded in $W_{0}$ and admits a convergent subsequence. We proceed by steps.\\

\noindent \textbf{Step 1.} The sequence $\{ u_{n}\}_{n \in \mathbb{N}}$ is bounded in $W_{0}$.\\
For $n$ large enough, by \ref{h1} in Definition \ref{PS}, we have
\begin{align*}
	E_{\lambda}^{p,1} (u_{n}) &= \Phi_{p,1}(u_{n}) - \lambda \Psi (u_{n})\\
	&= \dfrac{1}{p} \left( \displaystyle \int_{\mathbb{R}^{2N}}\dfrac{\vert u_{n}(x)-u_{n}(y) \vert^{p}}{\vert x-y \vert^{N+sp}}~dxdy + \int_{\Omega}\dfrac{\vert u_{n}(x)\vert^{p}}{\vert x \vert^{sp}}~dx\right) - \lambda \int_{\Omega}F(x,u_{n}(x))~dx\\
	&\leq d.
\end{align*}
In the other hand, by \ref{H2} we have
\begin{align*}
	E_{\lambda}^{p,1} (u_{n}) &\geq \dfrac{1}{p} \Vert u_{n}\Vert^{p}_{W_{0}} - \dfrac{\lambda}{\theta} \int_{\Omega}f(x,u_{n}(x))u_{n}(x)~dx\\
	&\geq \frac{1}{p} \Vert u_{n}\Vert^{p}_{W_{0}} - \frac{\lambda}{\theta} \int_{\Omega}f(x,u_{n}(x))u_{n}(x)~dx + \frac{1}{\theta} \Vert u_{n}\Vert^{p}_{W_{0}} - \frac{1}{\theta} \Vert u_{n}\Vert^{p}_{W_{0}}\\
	&> \left( \frac{1}{p} - \frac{1}{\theta}\right) \Vert u_{n}\Vert^{p}_{W_{0}} + \frac{1}{\theta} \left( \Vert u_{n}\Vert^{p}_{W_{0}} -\lambda \int_{\Omega}f(x,u_{n}(x))u_{n}(x)~dx\right)\\
	&\geq \left( \frac{1}{p} - \frac{1}{\theta}\right) \Vert u_{n}\Vert^{p}_{W_{0}} + \frac{1}{\theta} \langle \left(E_{\lambda}^{p,1}\right)^{'} (u_{n}), u_{n} \rangle.
\end{align*}
Due to \ref{h2} in Definition \ref{PS}, we have $\theta > p > 1$ and $\varepsilon \to 0$ such that
$$\Vert \left(E_{\lambda}^{p,1}\right)^{'} (u_{n})  \Vert_{W_{0}^{*}} \leq \varepsilon$$
and
\begin{align*}
	-\left| \frac{1}{\theta}\right| \langle \left(E_{\lambda}^{p,1}\right)^{'} (u_{n}), u_{n} \rangle  &\leq \Vert \left(E_{\lambda}^{p,1}\right)^{'} (u_{n})  \Vert_{W_{0}^{*}} \Vert u_{n}\Vert_{W_{0}}\\
	&\leq \varepsilon \Vert u_{n}\Vert_{W_{0}}.
\end{align*}
Thus,
$$\left( \frac{1}{p} - \frac{1}{\theta}\right) \Vert u\Vert^{p}_{W_{0}} \leq E_{\lambda}^{p,1} (u_{n}) - \frac{1}{\theta} \langle \left(E_{\lambda}^{p,1}\right)^{'} (u_{n}), u_{n} \rangle \leq d + \varepsilon \Vert u_{n}\Vert_{W_{0}}.$$
It  follows from this inequality that $\{ u_{n}\}$ is bounded in $W_{0}$.\\

\noindent \textbf{Step 2.} the sequence $\{ u_{n}\}_{n \in \mathbb{N}}$ possesses a convergent subsequence.\\
By the Eberlian-Smulyan theorem (see \cite[Theorem 21.D.]{Zeidler1990A}), passing to a subsequence if necessary, we can assume that $u_{n} \rightharpoonup u$. Then, because of the compactness of $\Psi^{'}$ thanks to condition \ref{H1} and to the compact embedding $W_{0} \hookrightarrow L^{q}(\Omega)$, for every $q \in [1,p^{*})$ we have
$$\Psi^{'}(u_{n}) \longrightarrow \Psi^{'}(u),$$
Since
$$\left(E_{\lambda}^{p,1}\right)^{'}(u_{n}) = \Phi_{p,1}^{'}(u_{n}) - \lambda \Psi^{'}(u_{n}) \longrightarrow 0,$$
then
$$\Phi_{p,1}^{'}(u_{n}) = \left(E_{\lambda}^{p,1}\right)^{'}(u_{n}) + \lambda \Psi^{'}(u_{n}) \longrightarrow 0 + \lambda \Psi^{'}(u).$$
Moreover, as $\Phi_{p,1}^{'}$ is a homeomorphism according to Lemma \ref{lemp}, then $u_{n} \rightarrow u$ in $W_{0}$ and so $E_{\lambda}^{p,1}$ satisfies the (PS)-condition.\hspace*{\fill} $\square$

	\begin{lem}
		\label{lem3}
		Assume that $f$ satisfies \ref{H2}. Then, there exist a positive $c$ constant such that
		\begin{equation}
			\label{estim}
			F(x,t) \geq c |t|^{\theta}, ~~~~\forall x \in \Omega, ~~|t|>M.
		\end{equation}
	\end{lem}
	
	\pr We start by setting $a(x) \coloneqq \displaystyle \min_{|\xi|= M} F(x, \xi)$ and
	\begin{equation}
		\label{varphi}
		\varphi_{t}(s) \coloneqq F(x,st), ~~~~\forall s>0.
	\end{equation}
	By the assumption \ref{H2}, we have for every $x \in \Omega$ and $|t|>M$ such that
	$$0 < \theta \varphi_{t}(s) = \theta F(x,st) \leq st f(x,st)=s\varphi_{t}^{'}(s), ~~\forall s > \frac{M}{|t|}.$$
	Therefore,
	$$\theta \varphi_{t}(s) \leq s\varphi_{t}^{'}(s).$$
	Also,
	$$\int^{1}_{\frac{M}{|t|}} \frac{\theta}{s}~ds \leq \int^{1}_{\frac{M}{|t|}} \frac{\varphi_{t}^{'}(s)}{\varphi_{t}(s)}~ds.$$
	We see that,
	$$\int^{1}_{\frac{M}{|t|}} \left( \ln |s|^{\theta}\right)^{'} ~ds \leq \int^{1}_{\frac{M}{|t|}} \left( \ln \left| \varphi_{t}(s)\right|\right)^{'} ~ds,$$
	
	$$\left[\ln |s|^{\theta} \right]^{1}_{\frac{M}{|t|}}  \leq \left[ \ln \left| \varphi_{t}(s)\right| \right]^{1}_{\frac{M}{|t|}},$$
	
	$$-\ln \left[\frac{M}{|t|} \right]^{\theta}   \leq  \ln \left| \varphi_{t}(1)\right| - \ln \left| \varphi_{t} \left( \frac{M}{|t|} \right)\right|.$$
	Then,
	$$\ln \varphi_{t} \left( \frac{M}{|t|} \right) - \ln \frac{M^{\theta}}{|t|^{\theta}}  \leq  \ln \varphi_{t}(1).$$
	We have
	$$- \ln \frac{M^{\theta}}{|t|^{\theta}} = -\left[ \ln M^{\theta} - \ln |t|^{\theta}\right] = \ln |t|^{\theta} - \ln M^{\theta} = \ln \frac{|t|^{\theta}}{M^{\theta}}.$$
	So,
	$$\ln \varphi_{t} \left( \frac{M}{|t|} \right) + \ln \frac{|t|^{\theta}}{M^{\theta}}  \leq  \ln \varphi_{t}(1),$$
	
	$$e^{\ln \varphi_{t} \left( \frac{M}{|t|} \right) + \ln \frac{|t|^{\theta}}{M^{\theta}}} = e^{\ln \left[ \varphi_{t} \left( \frac{M}{|t|} \right)  \frac{|t|^{\theta}}{M^{\theta}}\right]} \leq  e^{\ln \varphi_{t}(1)}.$$
	Thus,
	$$\varphi_{t} \left( \frac{M}{|t|} \right)  \frac{|t|^{\theta}}{M^{\theta}} \leq  \varphi_{t}(1).$$
	Taking into account of (\ref{varphi}), we obtain
	$$c |t|^{\theta}\leq a(x)\frac{|t|^{\theta}}{M^{\theta}}\leq F(x, \frac{M}{|t|} t)\frac{|t|^{\theta}}{M^{\theta}} \leq F(x,t),$$
	where $c\geq 0$ is a constant. Thus, (\ref{estim}) is proved.\hspace*{\fill} $\square$\\
	
	Now, we prove the following Proposition \ref{prop2} using the previous proved Lemma \ref{lem3}.
\begin{prop}
\label{prop2}
For every $\lambda > 0$ and the assumption \ref{H2} hold. Then, $E_{\lambda}^{p,1}$ is unbounded from below.
\end{prop}

\pr We fix $u_{0} \in W_{0}\backslash \{0\}$ and for each $t>1$, we have
\begin{align*}
	E_{\lambda}^{p,1}(tu_{0}) &= \Phi_{p,1} (tu_{0})-\lambda \Psi (tu_{0})\\
	&= \dfrac{t^{p}}{p} \left( \displaystyle \int_{\mathbb{R}^{2N}}\dfrac{\vert u_{0}(x)-u_{0}(y) \vert^{p}}{\vert x-y \vert^{N+sp}}~dxdy + \int_{\Omega}\dfrac{\vert u_{0}(x)\vert^{p}}{\vert x \vert^{sp}}~dx\right) - \lambda \int_{\Omega}F(x,u_{0}(x))~dx.
\end{align*}
So, from Lemma \ref{lem3} and Hardy's inequality (\ref{FHardy}) one has
\begin{align*}
	E_{\lambda}^{p,1}(tu_{0}) &\leq \dfrac{t^{p}(c_{H}+1)}{p} \Vert u_{0}\Vert^{p}_{W_{0}} - \lambda \int_{\Omega}c \left|tu_{0}(x)\right|^{\theta}~dx\\
	&\leq \dfrac{t^{p}(c_{H}+1)}{p} \Vert u_{0}\Vert^{p}_{W_{0}} - \lambda c t^{\theta} \displaystyle \int_{\Omega} \left|u_{0}(x)\right|^{\theta}~dx \\
	&< 0
\end{align*}
The assumption \ref{H2} ensures that $\theta>p$, this condition guarantees that
$$E_{\lambda}^{p,1}(tu_{0}) \longrightarrow -\infty, ~\text{ as } t \to +\infty,$$
which lead us to deduce that $E_{\lambda}^{p,1}$ is unbounded from below.\hspace*{\fill} $\square$\\
	
\prrr{\ref{thp1}} Setting $\lambda \in ] 0 ,\overline{\lambda} [$, we seek to apply Theorem \ref{th2} to problem (\ref{p+}) in the case $r=1$ to the space $X \coloneqq W_{0}$ and to the functionals
$$\Phi_{p,1} (u) \coloneqq \dfrac{1}{p} \left( \displaystyle \int_{\mathbb{R}^{2N}}\dfrac{\vert u(x)-u(y) \vert^{p}}{\vert x-y \vert^{N+sp}}~dxdy + \int_{\Omega}\dfrac{\vert u(x)\vert^{p}}{\vert x \vert^{sp}}~dx\right)$$
and
$$\Psi (u) \coloneqq \int_{\Omega}F(x,u(x))~dx.$$
The functional $\Phi_{p,1}$ is continuous and $\Phi^{'}_{p,1} : W_{0} \rightarrow W_{0}^{*}$ is a homeomorphism according to Proposition \ref{lemp}. Moreover, thanks to condition \ref{H1} and to the compact embedding $W_{0} \hookrightarrow L^{q}(\Omega)$, for every $q \in [1,p^{*})$. Then, $\Psi$ is continuous and has a compact derivative. \\

Via Proposition \ref{prop1} and Proposition \ref{prop2}, we established that $E_{\lambda}^{p,1}$ is unbounded from below and satisfies the Palais-Smale condition. Therefore, we can apply our main tool Theorem \ref{th2}, and it only remains to  prove that $\lambda \in ~]0,\overline{\lambda}[~\subseteq~]0,\frac{1}{\sup_{\{ \Phi_{p,1}(u)<1\}} \Psi(u)}[$ where $\overline{\lambda} \coloneqq \frac{q}{qa_{1}c_{1}p^{\frac{1}{p}} + a_{2}c^{q}_{q}p^{\frac{q}{p}}}$.

From (\ref{Heqp}), we have
$$\dfrac{1}{p} \Vert u\Vert^{p}_{W_{0}} \leq \Phi_{p,1}(u) < r, \text{  such that  } u \in \Phi_{p,1}^{-1}(]-\infty,r[).$$
For $r=1$, we get
$$\Vert u\Vert^{p}_{W_{0}} < p.$$
Then, for each $u \in W_{0}$
\begin{equation}
	\label{eq1}
	\Vert u\Vert_{W_{0}} < p^{\frac{1}{p}}, \text{  such that,  } u \in \Phi_{p,1}^{-1}(]-\infty,1[).
\end{equation}
Moreover, according to \ref{H1}, we have
\begin{align*}
	|F(x,t)| &= \left| \int_{0}^{t}f(x,s)~ds\right| \\
	&\leq \alpha \left| \int_{0}^{t}~ds\right| + \beta \left| \int_{0}^{t} |s|^{q-1}~ds\right|\\
	&=  \alpha |t| + \beta \dfrac{|t|^{q}}{q}.
\end{align*}
Then
\begin{align*}
	\Psi(u) &= \int_{\Omega} |F(x,u(x))|dx\\
	&\leq \int_{\Omega} \left[ \alpha |u(x)| + \beta \dfrac{|u(x)|^{q}}{q}\right]~dx \\
	&\leq \int_{\Omega} \alpha |u(x)|~dx + \dfrac{\beta}{q}\int_{\Omega} |u(x)|^{q}~dx \\	
	&= \alpha \Vert u\Vert_{L^{1}(\Omega)} + \dfrac{\beta}{q}  \Vert u\Vert_{L^{q}(\Omega)}^{q}
\end{align*}
Using the compact embedding $W_{0} \hookrightarrow L^{q}(\Omega)$ for every $q \in [1,p^{*})$ and for each $u \in \Phi_{p,1}^{-1}(]-\infty,1[)$, we have
$$\Psi(u) \leq  \alpha c_{1} \Vert u\Vert_{W_{0}} + \beta \left( c_{q} \Vert u\Vert_{W_{0}}\right)^{q}.$$
And by (\ref{eq1}), we obtain
\begin{align*}
\Psi(u) &<  \alpha c_{1} p^{\frac{1}{p}} + \dfrac{\beta}{q} c_{q}^{q} p^{\frac{q}{p}}\\
&= \dfrac{q \alpha c_{1} p^{\frac{1}{p}} + \beta c_{q}^{q} p^{\frac{q}{p}}}{q}.
\end{align*}
Since $\lambda \in ] 0 ,\overline{\lambda} [$, we get
$$\sup_{\{  u \in \Phi_{p,1}^{-1}(]-\infty,1[)\}} \Psi (u) < \dfrac{q \alpha c_{1} p^{\frac{1}{p}} + \beta c_{q}^{q} p^{\frac{q}{p}}}{q} \eqqcolon \dfrac{1}{\overline{\lambda}} < \dfrac{1}{\lambda},$$
from, the latter one has
$$0<\lambda<\overline{\lambda} \coloneqq \dfrac{q}{q \alpha c_{1} p^{\frac{1}{p}} + \beta c_{q}^{q} p^{\frac{q}{p}}} <\dfrac{1}{\sup_{\{ \Phi_{p,1}(u)<1\}} \Psi (u) }.$$
Then,
$$\lambda \in ]0,\overline{\lambda}[ \subseteq \left]0,\dfrac{1}{\sup_{\{ \Phi_{p,1}(u)<1\}} \Psi(u)}\right[.$$
Now that all hypotheses of Theorem \ref{th2} are verified. We conclude that for each $\lambda \in ]0,\overline{\lambda}[$, the functional $E_{\lambda}^{p,1}$ admits two distinct critical points that are weak solutions of problem (\ref{p+}). \hspace*{\fill} $\square$ 
\begin{rem}
Remark that Theorem \ref{thp1} confirm the existence of two positive weak solutions for problem (\ref{p+}), if the function $f$ is positive and $f(x,0) \neq 0$ in $\Omega$.
\end{rem}
The second result of this work is the following
\begin{theo}
	\label{thp2}
	Let $f : \Omega \times \mathbb{R} \rightarrow \mathbb{R}$ be a function with $f(x,0) \neq 0$ in $\Omega$ satisfying condition \ref{H1}. Then, there exists a positive number $\Lambda$ given by
$$\Lambda \coloneqq q\sup_{\{ \rho > 0\}} \left\{ \dfrac{\rho^{p-1}}{q\alpha c_{1} \left( \frac{p}{1-c_{H}}\right)^{\frac{1}{p}} + \beta c_{q}^{q} \left( \frac{p}{1-c_{H}}\right)^{\frac{q}{p}} \rho^{q-1}} \right\},$$
such that, for every $\lambda \in ] 0 , \Lambda [$, the following second case of  problem (\ref{P})
\leqnomode
\begin{equation*}
\label{p-}
\left\{\begin{array}{ll} (-\Delta )_p^s u - \dfrac{|u|^{p-2}u}{|x|^{sp}} = \lambda f(x,u) & \quad \mbox{in }\ \Omega\\[0.3cm]
			u= 0        & \quad \mbox{in }\ \mathbb{R}^N \setminus \Omega,\tag{P$_{-}$}
\end{array}
\right.
\end{equation*}
\reqnomode
admits at least one non-trivial weak solution $u_{\lambda} \in W_{0}$. Moreover,
$$\lim_{\lambda \rightarrow 0^{+}} \Vert u_{\lambda}\Vert_{W_{0}} = 0$$
and the function $g(\lambda) \coloneqq E_{\lambda}^{p,-1}(u_{\lambda})$ is negative and strictly decreasing in $] 0 , \Lambda [$.
\end{theo}
To achieve the proof of this result, we need to prove the following Lemmas:
\label{proof of lower semicontinuity}
\begin{lem}
\label{swls}
Let $s\in(0,1)$ and $N>ps$. Then, the functional $\Phi_{p,-1}$ is coercive and sequentially weakly lower semicontinious on $W_{0}$, i.e:
	$$\Phi_{p,-1} (u) \leq \liminf_{n \to +\infty} \Phi_{p,-1} (u_{n}) ~\text{ if }~ u_{n} \longrightarrow u ~\mbox{ weakly in }~W_{0}.$$
\end{lem}
\pr Using the fractional Hardy inequality (\ref{FHardy}), we obtain that for any $u \in W_{0}$,
\begin{align*}
	\Phi_{p,-1}(u) &= \dfrac{1}{p} \left( \displaystyle \int_{\mathbb{R}^{2N}}\dfrac{\vert u(x)-u(y) \vert^{p}}{\vert x-y \vert^{N+sp}}~dxdy - \int_{\Omega}\dfrac{\vert u(x)\vert^{p}}{\vert x \vert^{sp}}~dx\right)\\
	&\geq \left( \dfrac{1-c_{H}}{p}\right) \Vert u\Vert^{p}_{W_{0}}.
\end{align*}
We conclude that
$$\Phi_{p,-1}(u) \longrightarrow +\infty,\text{ as } \Vert u\Vert^{p}_{W_{0}} \longrightarrow +\infty,$$
which means that $\Phi_{p,-1} (u)$ is coercive on $W_{0}$.\\

Now, by \cite[Theorem 6.]{Fiscella2015} we know that $\mathbb{C}_c^{\infty}(\Omega)$ is dense subset of $W_{0}$. Hence, using density arguments, to prove that $\Phi_{p,-1}$ is sequentially weakly lower semicontinious on $W_{0}$, it is enough to show that the functional
\begin{equation}
\label{stat}
\Phi_{p,-1} \text{ is sequentially weakly lower semicontinious on } \mathbb{C}_c^{\infty}(\Omega).
	\end{equation}
	
	So, let $\{u_{n}\}_{n\in \mathbb{N}}$ be a sequence in $\mathbb{C}_c^{\infty}(\Omega)$ such that
\begin{equation}
\label{weakconv}
\begin{array}{ll}
u_{n} \longrightarrow u & \mbox{weakly in } ~W_{0}, \mbox{ as } n \rightarrow + \infty,
		\end{array}
\end{equation}
Thus, according to Theorem \ref{thccp} there exist two bounded measures $\mu$ and $\nu$, an at most enumerable set of indices $I$ of distinct points $\{x_{i}\}_{i\in I} \subset \overline{\Omega}$, and positive real numbers $\{ \mu_{i}\}_{i\in I}, ~\{ \nu_{i}\}_{i\in I}$, such that the following convergence hold weakly in the sense of measures,
\begin{equation}
	\label{w1}
\displaystyle \int_{\mathbb{R}^{N}}\dfrac{\vert u_{n}(x)-u_{n}(y) \vert^{p}}{\vert x-y \vert^{N+sp}}~dydx \rightharpoonup \mu \geq \displaystyle \int_{\mathbb{R}^{N}}\dfrac{\vert u(x)-u(y) \vert^{p}}{\vert x-y \vert^{N+sp}}~dydx + \sum_{i \in I} \mu_{i} \delta_{x_{i}},
\end{equation}
\begin{equation}
\label{w2}
\left|u_{n}(x)\right|^{p} \rightharpoonup \nu = \left|u(x)\right|^{p}dx + \sum_{i \in I} \nu_{i} \delta_{x_{i}},
\end{equation}
and finally,
\begin{equation}
\label{w3}
S_{p} ~ \nu_{i} \leq \mu_{i}, ~\forall i \in I,
\end{equation}
where  $S_{p}$ is the Sobolev constant given by (\ref{sobolev}) and $\delta_{x_{i}}$ denotes the Dirac mass at $x_{i}$.
	
From the continuity of the embedding $W_{0} \hookrightarrow L^{p}(\Omega)$, for every $p \in [1,p^{*}]$, we hve that
$$\begin{array}{ll}
u_{n} \longrightarrow u & \mbox{ stronly in }~L^{p}(\Omega), \mbox{ as } n \rightarrow + \infty.
\end{array} $$
By (\ref{w1}), (\ref{w2}) and (\ref{w3}) we get
\begin{align*}
\liminf_{n \to +\infty} \Phi_{p,-1}(u_{n}) &= \liminf_{n \to +\infty} \left[ \dfrac{1}{p} \displaystyle \int_{\mathbb{R}^{2N}}\dfrac{\vert u_{n}(x)-u_{n}(y) \vert^{p}}{\vert x-y \vert^{N+sp}}~dxdy - \dfrac{1}{p} \int_{\Omega}\dfrac{\vert u_{n}(x)\vert^{p}}{\vert x \vert^{sp}}~dx \right]\\
		&\geq \dfrac{1}{p} \left( \displaystyle \int_{\mathbb{R}^{2N}}\dfrac{\vert u(x)-u(y) \vert^{p}}{\vert x-y \vert^{N+sp}}~dxdy + \sum_{i=1}^{k} \mu_{i}\right) - \dfrac{1}{p} \left( \int_{\Omega}\dfrac{\vert u_{n}(x)\vert^{p}}{\vert x \vert^{sp}}~dx + \sum_{i=1}^{k} \nu_{i}\right)\\
		&= \dfrac{1}{p} \left( \displaystyle \int_{\mathbb{R}^{2N}}\dfrac{\vert u(x)-u(y) \vert^{p}}{\vert x-y \vert^{N+sp}}~dxdy - \int_{\Omega}\dfrac{\vert u(x)\vert^{p}}{\vert x \vert^{sp}}~dx\right) + \dfrac{1}{p} \sum_{i=1}^{k} \mu_{i} - \dfrac{1}{p} \sum_{i=1}^{k} \nu_{i}\\
		&\geq \Phi_{p,-1}(u) + \dfrac{1}{p} \sum_{i=1}^{k} \mu_{i} - \dfrac{1}{p} \sum_{i=1}^{k} \dfrac{\mu_{i}}{S_{p}} \\
		&= \Phi_{p,-1}(u) + \left( \dfrac{1-S_{p}^{-1} }{p}\right) \sum_{i=1}^{k} \mu_{i}\\
		&\geq \Phi_{p,-1}(u).
	\end{align*}
	This lead us to deduce statement stated in (\ref{stat}).
	
	Now, let $\{u_{n}\}_{n\in \mathbb{N}}$ be a sequence in $W_{0}$ satisfying the same condition (\ref{weakconv}). Then, using density arguments, we have for any $n\in \mathbb{N}$ there exists $\{u_{n}^{j}\}_{j \in \mathbb{N}} \in \mathbb{C}_{c}^{\infty}(\Omega)$ such that
\begin{equation}
		\label{weakconvj}
		\begin{array}{ll}
			u_{n}^{j} \longrightarrow u_{n} & \mbox{strongly in } ~W_{0}, \mbox{ as } j \rightarrow + \infty.
		\end{array}
	\end{equation}
From (\ref{weakconv}) and (\ref{weakconvj}), we have that for any $\varphi \in W_{0}$
$$\langle u_{n}^{j}-u,\varphi \rangle = \langle u_{n}^{j}-u_{n},\varphi\rangle + \langle u_{n}-u,\varphi\rangle \rightarrow 0, \mbox{ as } n,j \rightarrow + \infty.$$
Then,
	\begin{equation}
		\begin{array}{ll}
			u_{n}^{j} \longrightarrow u & \mbox{weakly in } ~W_{0}, \mbox{ as } n,j \rightarrow + \infty.
		\end{array}
	\end{equation}
	Since $\{u_{n}^{j}\}_{j \in \mathbb{N}} \in \mathbb{C}_{c}^{\infty}(\Omega)$ and the statment (\ref{stat}) is satisfied, we deduce that
	\begin{equation}
		\label{conv1}
		\liminf_{n,j \to +\infty} \Phi_{p,-1}(u_{n}^{j}) \geq  \Phi_{p,-1}(u).
	\end{equation}
	Moreover, by (\ref{weakconvj}) it is easy to see that for any $n \in \mathbb{N}$ we have
	$$\lim_{j \to +\infty} \Phi_{p,-1}(u_{n}^{j}) = \Phi_{p,-1}(u_{n}),$$
	so that, passing to $\liminf$ we get
	\begin{equation}
		\label{conv2}
		\liminf_{n,j \to +\infty} \Phi_{p,-1}(u_{n}^{j}) = \liminf_{n \to +\infty} \lim_{j \to +\infty} \Phi_{p,-1}(u_{n}^{j}) = \liminf_{n \to +\infty} \Phi_{p,-1}(u_{n})
	\end{equation}
	By (\ref{conv1}) and (\ref{conv2}) we get that
	$$\liminf_{n \to +\infty} \Phi_{p,-1}(u_{n}) \geq \Phi_{p,-1}(u)$$
	Therefore, $\Phi_{p,-1}$ is sequentially weakly lower semicontinious on $W_{0}$.\hspace*{\fill} $\square$
	
	\begin{lem}
		\label{swus}
		Let \ref{H1} be satisfied. Then, the functional $\Psi$ is sequentially weakly continuous.
	\end{lem}
	
	\pr let $\{u_{n} \}_{n \in \mathbb{N}}$ is bounded in $W_{0}$ and $W_{0}$ is a reflexive space. Then, up to a subsequence denoted by $u_{n}$, there exists $u \in W_{0}$, such that
	$$\begin{array}{ll}
		u_{n} \longrightarrow u & \mbox{ weakly in }~ W_{0},\\
		u_{n} \longrightarrow u & \mbox{ stronly in }~ L^{p}(\Omega),\\
		u_{n}(x) \longrightarrow u(x) & \mbox{ a.e. in }~\Omega,
	\end{array} $$
	Now, by assumption \ref{H1} we have that
	\begin{align*}
		|F(x,t)| &= \left| \int_{0}^{t}f(x,s)~ds\right| \nonumber \\
		&\leq \alpha \left| \int_{0}^{t}~ds\right| + \beta \left| \int_{0}^{t} |s|^{q-1}~ds\right| \nonumber\\
		&=  \alpha |t| + \beta \dfrac{|t|^{q}}{q}.
	\end{align*}
	further, according to the compact embedding $W_{0} \hookrightarrow L^{q}(\Omega)$ for every $q \in [1,p^{*})$, we get
	$$|\Psi (u_n)| \leq \int_{\Omega} |F(x,u_n)|~dx \leq  \alpha c_{1} \Vert u_n\Vert_{W_{0}} + \beta \left( c_{q} \Vert u_n\Vert_{W_{0}}\right)^{q} < + \infty,$$
	Thus, we apply the Lebesgue dominated convergence Theorem, we obtain
	$$\displaystyle \lim_{n \rightarrow 0} \int_{\Omega} \Psi (u_n)~dx = \int_{\Omega} \Psi (u)~dx.$$
	Then, the functional $\Psi$ is weakly semicontinuous.\hspace*{\fill} $\square$\\
	
\prrr{\ref{thp2}} We set $\lambda \in ] 0 , \Lambda [$. In order to apply Theorem \ref{th22} to problem (\ref{p-}) with the space $X = W_{0}$ and to the functionals
	$$\Phi_{p,-1} (u) \coloneqq \dfrac{1}{p} \left( \displaystyle \int_{\mathbb{R}^{2N}}\dfrac{\vert u(x)-u(y) \vert^{p}}{\vert x-y \vert^{N+sp}}~dxdy - \int_{\Omega}\dfrac{\vert u(x)\vert^{p}}{\vert x \vert^{sp}}~dx\right)$$
	and
	$$\Psi (u) \coloneqq \int_{\Omega}F(x,u(x))~dx.$$
	In view of Proposition \ref{lemp} and Lemma \ref{swls}, the functional $\Phi_{p,-1}$ is continuous, coercive and sequentially weakly lower semicontinious, also its $\displaystyle \inf_{u \in X}\Phi_{p,-1}(u)=0$. Moreover, the functional $\Psi$ is continuous, has a compact derivative and is sequentially weakly continuous according to Lemma \ref{swus}. We prove the theorem in the following steps.\\
	
	\label{stage1}
	\noindent \textbf{Step 1.} We start by proving that problem (\ref{p-}) admits at least one non-trivial weak solution $u_{\lambda} \in W_{0}$. \\
By \ref{H1}, we have
	\begin{align}
		\label{est1}
		|F(x,t)| &\leq \alpha |t| + \beta \dfrac{|t|^{q}}{q},~\forall (x,t) \in \Omega \times \mathbb{R}.
	\end{align}
Using the above inequality, we obtain
	\begin{align*}
		\Psi(u) &= \int_{\Omega} |F(x,u(x))|dx\\
		&\leq \int_{\Omega} \left[ \alpha |u(x)| + \beta \dfrac{|u(x)|^{q}}{q}\right]~dx \\
		&\leq \int_{\Omega} \alpha |u(x)|~dx + \dfrac{\beta}{q}\int_{\Omega} |u(x)|^{q}~dx \\	
		&= \alpha \Vert u\Vert_{L^{1}(\Omega)} + \dfrac{\beta}{q}  \Vert u\Vert_{L^{q}(\Omega)}^{q}.
\end{align*}
According to the compact embedding $W_{0} \hookrightarrow L^{q}(\Omega)$, for every $q \in [1,p^{*})$, we have
	\begin{align*}
		\Psi(u) &\leq  \alpha c_{1} \Vert u\Vert_{W_{0}} + \beta \left( c_{q} \Vert u\Vert_{W_{0}}\right)^{q}.
	\end{align*}
On the other hand, we get from (\ref{Heqm}) that
	\begin{equation}
		\label{s2}
		\Vert u\Vert_{W_{0}} < \left( \dfrac{pr}{1-c_{H}}\right)^{\frac{1}{p}},~\forall u \in W_{0},~\Phi_{p,-1}(u)<r.
	\end{equation}
Now, from (\ref{s2}), one has
\begin{align*}
\Psi(u) &< \alpha c_{1} \left( \dfrac{pr}{1-c_{H}}\right)^{\frac{1}{p}} + \beta \dfrac{c_{q}^{q}}{q} \left( \dfrac{pr}{1-c_{H}}\right)^{\frac{q}{p}},
\end{align*}
for every $u \in W_{0}$ such that $\Phi_{p,-1}(u)<r$.

\noindent Then,
$$\sup_{\{  u \in \Phi_{p,-1}^{-1}(]-\infty,r[)\}}\Psi(u) < \alpha c_{1} \left( \dfrac{p}{1-c_{H}}\right)^{\frac{1}{p}} r^{\frac{1}{p}} + \beta \dfrac{c_{q}^{q}}{p} \left( \dfrac{p}{1-c_{H}}\right)^{\frac{q}{p}} r^{\frac{q}{q}}.$$
Hence, for every $r \in ]0,+\infty[$ we have
$$\dfrac{\left( \sup_{\{ u \in \Phi_{p,-1}^{-1}( ]-\infty,r[ )\}} \Psi(u) \right) }{r} \leq  \alpha c_{1} \left( \dfrac{p}{1-c_{H}}\right)^{\frac{1}{p}}r^{\frac{1-p}{p}} + \beta \dfrac{c_{q}^{q}}{q} \left( \dfrac{p}{1-c_{H}}\right)^{\frac{q}{p}} r^{\frac{q-p}{p}}.$$
In particular, for $r=\rho^{p}$ we have
\begin{equation}
\label{s3}
\dfrac{\left( \sup_{\{ u \in \Phi_{p,-1}^{-1}( ]-\infty,\rho^{p}[ )\}} \Psi(u) \right) }{\rho^{p}}\leq \alpha c_{1} \left( \dfrac{p}{1-c_{H}}\right)^{\frac{1}{p}} \rho^{1-p} + \beta \dfrac{c_{q}^{q}}{q} \left( \dfrac{p}{1-c_{H}}\right)^{\frac{q}{p}} \rho^{q-p}.
\end{equation}
Now, setting the function $u_{0} \in W_{0}$ such that
\begin{equation}
\label{z}
 u_{0} \in \Phi_{p,-1}^{-1}( ]-\infty,\gamma^{p}[ ).
\end{equation}
Also, we observe that
\begin{equation}
\label{z1}
\Phi_{p,-1}(u_{0}) = \dfrac{1}{p} \Vert u_{0}\Vert^{p}_{W_{0}} + \dfrac{1}{p} \int_{\Omega}\dfrac{\vert u_{0}(x)\vert^{p}}{\vert x \vert^{sp}}~dx = 0
\end{equation}
and
\begin{equation}
\label{z2}
\Psi(u_{0}) = \int_{\Omega}F(x,u_{0}(x))~dx = 0.
\end{equation}
Therefore, from (\ref{z}), (\ref{z1}) and (\ref{z2}), we obtain that
\begin{align*}
\varphi(\rho^{p}) &\coloneqq \displaystyle \inf_{u \in \Phi_{p,-1}^{-1}( ]-\infty,\rho^{p}[ )} \dfrac{\left( \sup_{\{ v \in \Phi_{p,-1}^{-1}( ]-\infty,\rho^{p}[ )\}} \Psi(v) \right)-\Psi(u)}{\rho^{p}-\Phi_{p,-1}(u)}\\
&\leq \dfrac{\left( \sup_{\{ v \in \Phi_{p,-1}^{-1}( ]-\infty,\rho^{p}[ )\}} \Psi(v) \right)-\Psi(u_{0})}{\rho^{p}-\Phi_{p,-1}(u_{0})}\\
&=\dfrac{\left( \sup_{\{ v \in \Phi_{p,-1}^{-1}( ]-\infty,\rho^{p}[ )\}} \Psi(v) \right)}{\rho^{p}}.
\end{align*}
Thus, by (\ref{s3})we have
\begin{align*}
\varphi(\rho^{p}) &\leq \alpha c_{1} \left( \dfrac{p}{1-c_{H}}\right)^{\frac{1}{p}} \rho^{1-p} + \beta \dfrac{c_{q}^{q}}{q} \left( \dfrac{p}{1-c_{H}}\right)^{\frac{q}{p}} \rho^{q-p}.
\end{align*}
Furthermore, since $0<\lambda <\Lambda$
$$\varphi(\rho^{p}) \leq  \alpha c_{1} \left( \dfrac{p}{1-c_{H}}\right)^{\frac{1}{p}} \rho^{1-p} + \beta \dfrac{c_{q}^{q}}{q} \left( \dfrac{p}{1-c_{H}}\right)^{\frac{q}{p}} \rho^{q-p} \eqqcolon \dfrac{1}{\Lambda (\rho)} < \dfrac{1}{\lambda}.$$
So,
$$0<\lambda<\Lambda (\rho) \coloneqq \dfrac{q\rho^{p-1}}{q\alpha c_{1} \left( \frac{p}{1-c_{H}}\right)^{\frac{1}{p}} + \beta c_{q}^{q} \left( \frac{p}{1-c_{H}}\right)^{\frac{q}{p}} \rho^{q-1}} \leq \dfrac{1}{\varphi(\rho^{p})}.$$
Then,
$$\lambda \in \left]0,\Lambda (\rho) \right[ \subseteq \left]1,\dfrac{1}{\varphi(\rho^{p})} \right[$$
	In conclusion, according to Theorem \ref{th22}, there exists a critical point $u_{\lambda} \in \Phi_{p,-1}^{-1}( ]-\infty,\rho^{p}[ )$ for $E_{\lambda}^{p,-1}$ in $W_{0}$ which is a global minimum of the restriction $E_{\lambda}^{p,-1}$ to $\Phi_{p,-1}^{-1}( ]-\infty,\rho^{p}[ )$. Moreover, the function $u_{\lambda}\neq0$ since $f(x,0)\neq 0$ in $\Omega$.\\
	
\label{stage2}
\noindent \textbf{Step 2.} We show that $\lim_{\lambda \to 0^{+}}\Vert u_{\lambda}\Vert_{W_{0}} = 0$ and that the functional $g(\lambda) \coloneqq E_{\lambda}^{p,-1}(u_{\lambda})$ is negative and strictly decreasing in $\left]0,\Lambda(\rho) \right[$.\\
As $\Phi_{p,-1}$ is coercive, then $u_{\lambda} \in \Phi_{p,-1}^{-1} ( ]-\infty,\rho^{p}[ )$ is bounded in $W_{0}$, that is to say
$$\Vert u_{\lambda} \Vert_{W_{0}} \leq K, ~\text{for } K>0$$
Thus, due to the compactness of the operator $\Psi^{'}$, there exists a constant $C>0$ such that
\begin{equation}
\label{est2}
\left| \langle \Psi^{'} (u_{\lambda}),u_{\lambda}\rangle\right| \leq \Vert \Psi^{'} (u_{\lambda}) \Vert_{W_{0}^{*}} \Vert u_{\lambda} \Vert_{W_{0}}<CK^{2}, ~\forall \lambda \in \left]0,\Lambda(\rho)\right[.
\end{equation}
On the other hand, since $u_{\lambda}$ is a critical point of $E_{\lambda}^{p,-1}$, then
	$$\langle \left( E_{\lambda}^{p,-1}\right)^{'} (u_{\lambda}),u_{\lambda}\rangle = 0,$$
which implies that
\begin{equation*}
\langle \Phi_{p,-1}^{'} (u_{\lambda}) - \lambda \Psi^{'} (u_{\lambda}),u_{\lambda}\rangle = 0, ~\forall \lambda \in \left]0,\Lambda(\rho) \right[.
\end{equation*}
So,
\begin{equation}
\label{est3}
p \Phi_{p,-1}(u_{\lambda})=\langle \Phi_{p,-1}^{'} (u_{\lambda}),u_{\lambda}\rangle = \lambda \langle \Psi^{'} (u_{\lambda}),u_{\lambda}\rangle, ~\forall \lambda \in \left]0,\Lambda(\rho) \right[.
\end{equation}
Therefore, by (\ref{est2}) and (\ref{est3}), we get
\begin{equation}
\label{est4}
\lim_{\lambda \rightarrow 0^{+}} p\Phi_{p,-1} (u_{\lambda}) = \lim_{\lambda \rightarrow 0^{+}} \lambda \langle \Psi^{'} (u_{\lambda}),u_{\lambda}\rangle = 0, ~\forall p<1.
\end{equation}
Moreover, by (\ref{Heqm}) one has
\begin{equation}
\label{est5}
\Vert u_{\lambda} \Vert^{p}_{W_{0}} \leq  \dfrac{p\Phi_{p,-1}(u_{\lambda})}{1-c_{H}},~\forall \lambda \in \left]0,\Lambda(\rho) \right[.
	\end{equation}
Then, we conclude by the conditions (\ref{est4}) and (\ref{est5}) that
$$\lim_{\lambda \rightarrow 0^{+}} \Vert u_{\lambda}\Vert_{W_{0}} = 0$$
Furthermore, since the restriction $E_{\lambda}^{p,-1}$ to $\Phi_{p,-1}^{-1}( ]-\infty,\gamma^{p}[ )$ admits a global minimum, which is a local minimum of $E_{\lambda}^{p,-1}$ in $W_{0}$, the map $g(\lambda) \coloneqq E_{\lambda}^{p,-1}(u_{\lambda})$ is negative in $\left] 0 , \Lambda(\rho) \right[$, because $u_{\lambda} \neq 0$ and $E_{\lambda}^{p,-1}(0)=0$.
	
	Finally, we prove that the function $g(\lambda) \coloneqq E_{\lambda}^{p,-1}(u_{\lambda})$ is negative and strictly decreasing in $] 0 , \Lambda [$, by observing that
$$E_{\lambda}^{p,-1} (u) =\lambda \left( \dfrac{\Phi_{p,-1} (u)}{\lambda} -  \Psi (u)\right).$$
Now, we assume $u_{\lambda_1},u_{\lambda_2} \in W_{0}$ are critical points of $E_{\lambda}^{p,-1}$, for every $\lambda_1,\lambda_2 \in \left] 0 , \Lambda(\rho) \right[$, with $\lambda_1<\lambda_2$. Further, we set
$$I_{\lambda_{i}} \coloneqq \inf_{u \in \Phi_{p,-1}^{-1}( ]-\infty,\gamma^{p}[ )} \left( \dfrac{\Phi_{p,-1} (u)}{\lambda_{i}} -  \Psi (u)\right)=\dfrac{1}{\lambda_{i}} E_{\lambda_{i}}^{p,-1} (u_{\lambda_{i}}), ~i=1,2.$$
Obviously, as mentioned earlier $I_{\lambda_{i}}<0$ for $i=1,2$, and since $\lambda_1<\lambda_2$, we have $I_{\lambda_{2}}\leq I_{\lambda_{1}}$. Therefore,
	$$E_{\lambda_{2}}^{p,-1} (u_{\lambda_{2}})=\lambda_{2} I_{\lambda_{2}} \leq \lambda_{2} I_{\lambda_{1}} < \lambda_{1} I_{\lambda_{1}} = E_{\lambda_{1}}^{p,-1} (u_{\lambda_{1}}).$$
	We conclude that, as $\lambda \in \left]0,\Lambda \right[$ is arbitrary, the above conclusions are still true in $\left]0,\Lambda \right[$. Hence the proof is completed.\hspace*{\fill} $\square$
	
	\begin{rem}
	To calculate the maximum of $\Lambda$, we need to look at its first derivative
	$$\Lambda^{'}(\rho) = - \left[ \dfrac{ \alpha c_{1} \left( \dfrac{p}{1-c_{H}}\right)^{\frac{1}{p}} (1-p)\rho^{-p} + \beta \dfrac{c_{q}^{q}}{q} \left( \dfrac{p}{1-c_{H}}\right)^{\frac{q}{p}} (q-p) \rho^{q-p-1} }{\left( \alpha c_{1} \left( \dfrac{p}{1-c_{H}}\right)^{\frac{1}{p}} \rho^{1-p} + \beta \dfrac{c_{q}^{q}}{q} \left( \dfrac{p}{1-c_{H}}\right)^{\frac{q}{p}} \rho^{q-p}\right)^{2}}\right].$$
	Further, we set $\Lambda^{'}(\rho)$ equal to zero and obtain
	$$\alpha c_{1} \left( \dfrac{p}{1-c_{H}}\right)^{\frac{1}{p}} (1-p)\rho^{-p} + \beta \dfrac{c_{q}^{q}}{q} \left( \dfrac{p}{1-c_{H}}\right)^{\frac{q}{p}} (q-p) \rho^{q-p-1}=0.$$
	Then,
	$$\rho_{\max} \coloneqq \left( \dfrac{p}{1-c_{H}}\right)^{-\frac{1}{p}} \left[ q \dfrac{\alpha c_{1}}{\beta c_{q}^{q}} \left( \dfrac{1-p}{p-q}\right)  \right]^{\frac{1}{q-1}}. $$
	Which lead us to conclude that $\Lambda$ is defined as follows
	\begin{equation*}
	 	\Lambda(\rho) = \left\{ \begin{array}{ll}
	 		+\infty & \quad \mbox{if }\  1<q<p \\
			\frac{1-c_{H}}{\beta c_{q}^{q}}      & \quad \mbox{if }\ q=p \\
			\frac{q\rho_{\max}^{p-1}}{q\alpha c_{1} \left( \frac{p}{1-c_{H}}\right)^{\frac{1}{p}} + \beta c_{q}^{q} \left( \frac{p}{1-c_{H}}\right)^{\frac{q}{p}} \rho_{\max}^{q-1}}  & \quad \mbox{if }\ q \in ]p,p^{*}[.
		\end{array}
		\right.
	\end{equation*}
	Notice from the latter that if $f$ satisfy the condition \ref{H1} at infinity, with $q \in ]1,p[$. Then, from Theorem \ref{thp2} we confirm that for each $\lambda>0$, our problem (\ref{P-}) admits at least a non-trivial weak solution.
\end{rem}

\end{document}